\documentclass{amsart}

\usepackage[utf8]{inputenc}

\usepackage{url}

\usepackage[all]{xy}
\usepackage{pb-diagram, pb-xy}

\usepackage{amssymb}
\usepackage{amsthm}

\newtheorem{thm}{Theorem}[section]
\newtheorem{lem}[thm]{Lemma}
\newtheorem{cor}[thm]{Corollary}
\newtheorem{prop}[thm]{Proposition}

\newcommand{\bbC}{{\mathbb C}}

\newcommand{\bbH}{{\mathbb H}}
\newcommand{\bbL}{{\mathbb L}}
\newcommand{\bbN}{{\mathbb N}}

\newcommand{\bbQ}{{\mathbb Q}}
\newcommand{\bbR}{{\mathbb R}}

\newcommand{\bbV}{{\mathbb V}}
\newcommand{\bbZ}{{\mathbb Z}}

\newcommand{\calA}{{\mathcal A}}
\newcommand{\calC}{{\mathcal C}}
\newcommand{\calH}{{\mathcal H}}

\newcommand{\calM}{{\mathcal M}}
\newcommand{\calO}{{\mathcal O}}
\newcommand{\calR}{{\mathcal R}}
\newcommand{\calT}{{\mathcal T}}
\newcommand{\calX}{{\mathcal X}}

\newcommand{\frakA}{{\mathfrak A}}
\newcommand{\frakS}{{\mathfrak S}}

\newcommand{\Dbc}{{{D^{\hspace{0.01em}b}_{\hspace{-0.13em}c} \hspace{-0.05em} }}}
\newcommand{\pH}{{^p \! H}}
\newcommand{\pdelta}{{^p\hspace{-0.05em} \delta}}

\newcommand{\Perv}{\mathrm{Perv}}
\newcommand{\IC}{{\mathrm{IC}}}

\newcommand{\Sp}{\mathit{Sp}}

\newcommand{\Sl}{\mathit{Sl}}
\newcommand{\Pic}{{\mathit{Pic}}}
\newcommand{\Hom}{{\mathit{Hom}}}
\newcommand{\Aut}{{\mathit{Aut}}}
\newcommand{\id}{{\mathit{id}}}
\newcommand{\diag}{\mathit{diag}}
\newcommand{\sgn}{{\mathit{sgn}}}

\begin{document}

\title[On a family of surfaces of general type]{On a family of surfaces of general type attached to abelian fourfolds
 and the Weyl group~$W(E_6)$} 
\author{T. Kr\"amer}
\address{Mathematisches Institut, Ruprecht-Karls-Universit\"at Heidelberg\\ Im Neuenheimer Feld 288, D-69120 Heidelberg (Germany)}
\email{tkraemer@mathi.uni-heidelberg.de}


\begin{abstract}
We study a family of surfaces of general type that arises from the intersections of two translates of the theta divisor on a principally polarized complex abelian fourfold. In particular we determine the N\'eron-Severi lattices of these surfaces and show that for a general abelian fourfold the monodromy group of the associated variation of Hodge structures is a subgroup of index at most two in the Weyl group $W(E_6)$, whereas for a Jacobian variety it degenerates to a subgroup of the symmetric group $\frakS_6$.
\end{abstract}

\maketitle

\thispagestyle{empty}

\section{Introduction}

Let $X$ be a principally polarized abelian variety (ppav) of dimension~$g>1$ over the complex numbers. The polarization determines an ample divisor up to a translation by a point $x\in X(\bbC)$, and we choose the theta divisor $\Theta \subset X$ to be one of the~$2^{2g}$ translates which are symmetric in the sense that they are stable under the inversion morphism $-\id_X: X\rightarrow X$. For any point $x$ we denote by~$\Theta_x = \Theta + x$ the corresponding translate of our theta divisor. The geometry of the intersections
\[
 Y_{x} \;=\; \Theta \cap \Theta_{x} 
\]
is closely connected with the moduli of the ppav and has been studied in relation with Torelli's theorem~\cite{Weil} \cite{DebW}, with the Schottky problem~\cite{BeNovikov} \cite{DeUpdate} \cite{GruSchottky} and with the Prym map~\cite{Do} \cite{Iz}. One of our goals in this paper is to see how the cohomology of these intersections varies with the point~$x$, with a focus on the case $g=4$ which leads to a family of surfaces of general type with a rich geometric structure.

\medskip

In section~\ref{sec:intersections} we discuss some general features of the above intersections. We always have an involution $\sigma_x:  Y_x \rightarrow Y_x, \, y \mapsto  x-y$ by the symmetry of the theta divisor. If the theta divisor is smooth or has at most isolated singularities, then for general $x$ a Bertini-type argument shows that $Y_x$ is smooth of dimension $g-2$ and that~$\sigma_x$ is \'etale so that the quotient
\[ Y_x^+ \;=\; Y_x / \sigma_x \] 
is again a smooth variety of the same dimension. These constructions are interesting already for $g=3$ where they show that the fibres of the Prym map over $\calA_3$ are birational to Kummer varieties~\cite{Re} \cite[\S 6.4]{Do}. After a brief review of this case, we focus on the case $g=4$. Here the Hirzebruch-Riemann-Roch theorem generically gives the following numerical data. \pagebreak

\begin{lem}
Let $X$ be a complex ppav of dimension $g=4$ with a smooth theta divisor, and fix a general point $x\in X(\bbC)$. Then $Y_x$ and $Y_x^+$ are smooth surfaces of general type with the following Hodge diamonds. \medskip
\[
h^{p,q}(Y_x) \quad
\xymatrix@R=0.05em@C=0.05em{
 && 1 && \\
 & 4 && 4 & \\
 17 && 52 && 17 \\
 & 4 && 4 & \\
 && 1 && 
}
\qquad \qquad \qquad
h^{p,q}(Y_x^+) \quad
\xymatrix@R=0.05em@C=0.05em{
 && 1 && \\
 & 0 && 0 & \\
 6 && 22 && 6 \\
 & 0 && 0 & \\
 && 1 && 
}
\]
\end{lem}

\noindent
The details of the computation can be found in lemma~\ref{lem:hodge-numbers}. Since $h^{1,1}(X)=16$ and $h^{2,0}(X)=h^{0,2}(X)=6$, the Lefschetz hyperplane theorem then implies that in dependence of the point $x$ the varying part of the cohomology $H^2(Y_x^+, \bbZ)$ defines over some Zariski-open dense subset $U\subset X$ a variation of Hodge structures of pure Hodge type $(1,1)$ and rank six. This variation of Hodge structures will be one of our main objects of study, and we denote it by $\bbV_+$ in what follows.

\medskip

With the above information the Hodge index theorem shows that the intersection form has signature $(35,51)$ on $H^2(Y_x, \bbQ)$ and signature~$(13, 21)$ on $H^2(Y_x^+, \bbQ)$, see corollary~\ref{cor:signatures}. Via a spectral sequence argument we also determine the integral cohomology groups of our surfaces (including torsion) in proposition~\ref{prop:coh}. As usual we consider the middle cohomology groups modulo torsion as integral lattices with respect to the intersection form. In terms of the standard hyperbolic plane $U$ these lattices have the following shape as we will see in section~\ref{sec:lattice}.

\begin{prop} 
If $X$ is a complex ppav of dimension $g=4$ with a smooth theta divisor, then for a general point $x\in X(\bbC)$ we have
 \begin{eqnarray} \nonumber
 H^2(Y_x, \bbZ) &\cong &  U^{35} \oplus E_8^2(-1), \\ \nonumber
 H^2(Y_x^+, \bbZ)/\mathit{Torsion} &\cong &  U^{13} \oplus E_8(-1). 
 \end{eqnarray}
\end{prop}

\noindent
A closer look at the configuration of these lattices in proposition~\ref{prop:E6} will show that the varying part of $H^2(Y_x^+, \bbZ)$ is a rescaled copy of the $E_6$ lattice. This being said, let us consider the monodromy group $G$ of the local system underlying the variation~$\bbV_+$ of Hodge structures over the Zariski-open dense subset $U\subset X$ from above. By definition this monodromy group is the image of the representation of the fundamental group~$\pi_1(U, x)$ on the stalk of~$\bbV_+$ at a chosen point $x\in U(\bbC)$. Up to isomorphism this image is of course independent of our choices, so we suppress the point $x$ from the notation. Since for any smooth proper family of varieties the monodromy operation preserves the intersection form on the fibres, our lattice computations imply that $G$ must be a subgroup of the Weyl group $W(E_6)$, see proposition~\ref{prop:upper-bound}. One of the main results of this paper is the following statement which goes back to a conjecture of R. Weissauer.

\begin{thm} \label{thm:generic_monodromy} 
For a general complex ppav $X$ of dimension $g=4$, the monodromy group $G$ of the local system underlying $\bbV_+$ is either the Weyl group $W(E_6)$ or its unique simple subgroup of index two which is the kernel of the sign homomorphism $\sgn: W(E_6) \longrightarrow \{ \pm 1\}$.
\end{thm}

We remark that the occurance of the Weyl group $W(E_6)$ in the above context is related to the~$27$ lines on a cubic surface~\cite{Do}. The link is given by the work of E.~Izadi on Prym-embedded curves in abelian fourfolds~\cite{Iz}. We discuss this relationship in section~\ref{sec:upper-bound}, though it will not be used explicitly in the sequel. 

\medskip

Let us briefly mention the main ideas in the proof of theorem~\ref{thm:generic_monodromy}, referring to section~\ref{sec:translates} for more details. We have already observed that the monodromy group~$G$ must be a subgroup of the Weyl group $W(E_6)$, so all we have to do is to show that for a general ppav $X$ this monodromy group is sufficiently large. For this we use two ideas. On the one hand, in section~\ref{sec:irreducible} we relate the variation of Hodge structures~$\bbV_+$ to the convolution square of the theta divisor in the sense of~\cite{KrWVanishing}. By the Tannakian formalism of loc.~cit.~the decomposition of such convolutions into irreducible pieces is controlled by the representation theory of a certain algebraic group attached to the theta divisor. For a general ppav this group has been determined in~\cite{KrWSchottky}, and this will imply that the local system underlying~$\bbV_+$ is irreducible so that the subgroup~$G \subseteq W(E_6)$ must act irreducibly on the vector space~$E_6 \otimes_\bbZ \bbC$. On the other hand, we study a degeneration of our general ppav~$X$ into the Jacobian variety~$JC$ of a smooth projective curve~$C$. For Jacobian varieties we have the following result that may be of interest in its own right.

\medskip

Let $C$ be a smooth complex projective curve of even genus $g=2n$, and identify the symmetric product
\[ C_n \;=\; (C\times \cdots \times C)/\frakS_n  \]
with the set of effective divisors of degree $n$ on the curve. We can then consider the difference morphism
\[
 d_n: \quad C_n \times C_n \; \longrightarrow \; JC \;=\; \Pic^0(C)
\]
which sends a pair $(D, E)$ of effective divisors of degree $n$ on the curve $C$ to the isomorphism class of the line bundle $\calO_C(D-E)$. A simple calculation of intersection numbers shows that $d_n$ is generically finite of degree $N=\tbinom{2n}{n}$. Thus over some Zariski-open dense subset of $JC$ the difference morphism restricts to a finite \'etale cover, and we can consider its Galois group $G(d_n)$. For a general curve $C$ this group will be determined in section~\ref{sec:diff} via a degeneration into a hyperelliptic curve, with the following overall result.

\begin{thm} \label{thm:general-difference}
For a general curve $C$ of even genus $g=2n$, the Galois group~$G(d_n)$ is either the alternating group $\frakA_N$ or the full symmetric group $\frakS_N$.
\end{thm}

Returning again to the monodromy group $G$ attached to a general ppav $X$, a degeneration into a general Jacobian variety will show that for $g=4$ we have an embedding $G(d_2)\hookrightarrow G$. This lower bound will be sufficient to conclude the proof of theorem~\ref{thm:generic_monodromy}. We also remark that if in this argument the Galois group happens to be the full symmetric group~$\frakS_6$, then the monodromy group $G$ has to be the full Weyl group~$W(E_6)$ because it then contains a reflection. However, at present we do not know whether this is the case generically.

\medskip

\section{Intersections of theta divisors} \label{sec:intersections}

Let $X$ be a complex ppav of dimension $g\geq 2$. Throughout this section we will always assume that our chosen symmetric theta divisor $\Theta \subset X$ is smooth or has at most isolated singularities so that we have the following Bertini-type result.

\begin{lem} \label{lem:bertini}
If the theta divisor has at most isolated singularities, then there exists a Zariski-open dense subset $U\subset X$ and a smooth proper family
\[ f_U: \;  Y_U \;  \longrightarrow \; U \]
whose fibre over any $x\in U(\bbC)$ is isomorphic to $Y_x = \Theta \cap \Theta_x$. 
\end{lem}

{\em Proof.} Recall from~\cite[cor.~III.10.7]{Hartshorne} that if $f: V \rightarrow W$ is a morphism of complex algebraic varieties and if $V$ is smooth, then over some open dense $U\subseteq W$ the restriction
$
 f_U:  f^{-1}(U) \rightarrow U
$
is a smooth morphism. Shrinking $W$ we can in fact replace the smoothness condition on $V$ by the weaker condition that the singular locus $\mathit{Sing}(V)$ is mapped via~$f$ into a proper closed subset of $W$. We apply this to the addition morphism 
\[
f: \quad V \;=\; \Theta \times \Theta \; \longrightarrow \; X \;=\; W, \quad (y,z) \mapsto y+z.
\]
Here $\mathit{Sing}(V) = (\mathit{Sing}(\Theta) \times \Theta) \cup (\Theta \times \mathit{Sing}(\Theta))$ has dimension~$g-1$ or is empty.~In both cases we obtain that over some open dense $U\subset X$ the restriction~$f_U$ is smooth of relative dimension $g-2$. By construction we have $f^{-1}(x) \cong \Theta \cap \Theta_{x}$ for all $x\in X(\bbC)$, hence the claim of the lemma follows. \qed 

\medskip

The symmetry of the theta divisor implies that the intersections $Y_x = \Theta \cap \Theta_x$ are stable under the involution $ \sigma_x: X \rightarrow X, y  \mapsto x-y$, and we obtain

\begin{cor} \label{cor:smooth-and-etale}
If the theta divisor has at most isolated singularities, then there is a Zariski-open dense $U\subset X$ such that for all $x\in U(\bbC)$ the quotient morphism
\[
   Y_x \; \longrightarrow \; Y_x^+ \; = \; Y_x/\sigma_x
\]
is an \'etale double covering between smooth varieties of dimension $g-2$.
\end{cor}

{\em Proof.} Lemma~\ref{lem:bertini} gives an open dense subset $U\subset X$ such that $Y_x$ is smooth of dimension $g-2$ for all $x\in U(\bbC)$. Shrinking this open dense subset we can assume that it does not contain points of the form $x=2y$ with $y\in \Theta(\bbC)$. Then the involution $\sigma_x: Y_x \rightarrow Y_x$ is \'etale for all~$x\in U(\bbC)$. 
\qed

\medskip

For the rest of this section we always fix a point~$x$ with the properties in the above corollary, and we put 
\[ Y \;=\;Y_x, \quad Y^+ \;=\; Y_x^+  \quad \textnormal{and} \quad \sigma \;=\; \sigma_x \]
for brevity. Like for any \'etale double cover of smooth complex varieties, the rational cohomology of the quotient $Y^+$ can be naturally identified with the eigenspace for the eigenvalue $+1$ in
\[
 H^\bullet(Y, \bbQ) \;=\; H^\bullet(Y, \bbQ)^+ \oplus H^\bullet(Y, \bbQ)^-
\]
where the upper index $\, \pm \,$ indicates that $\sigma$ acts by the scalar~$\pm 1$ on the respective two eigenspaces. Now a version of the weak Lefschetz theorem~\cite[rem.~3.1.29]{LazarsfeldPositivityI} says that if $W$ is a smooth complex projective variety of dimension $d$, then for all ample effective divisors $D_1, \dots, D_r$ the restriction map $H^n(W, \bbZ) \longrightarrow H^n(D_1\cap \cdots \cap D_r, \bbZ)$ is an isomorphism for~$n<d-r$ and a monomorphism for $n=d-r$. Note that the intersecting divisors are not required to be smooth or transverse to each other. In our case it follows that the restriction morphism
\[
 H^n(X, \bbQ) \; \longrightarrow \; H^n(Y, \bbQ)
\]
is an isomorphism for $n<g-2$ and a monomorphism for $n=g-2$. So by Poincar\'e duality the interesting part of the cohomology of $Y$ sits in degree $n=g-2$ and can be defined as the orthocomplement
\[
 V \;=\; H^{g-2}(X, \bbQ)^\perp \;\subseteq\; H^{g-2}(Y, \bbQ)
\]
with respect to the intersection form. Now the involution $\sigma$ acts by $\epsilon = (-1)^{g}$ on~$H^{g-2}(X, \bbQ)$, hence
\[
 H^{g-2}(X, \bbQ) \; \subseteq \; H^{g-2}(Y, \bbQ)^\epsilon \;=\;
 \begin{cases}
 H^{g-2}(Y, \bbQ)^+ & \textnormal{for $g$ even}, \\
 H^{g-2}(Y, \bbQ)^- & \textnormal{for $g$ odd}.
 \end{cases}
\]
So the orthocomplement from above admits a decomposition $V=V_+\oplus V_-$ into the eigenspaces
\[
 V_{-\epsilon} \;=\; H^{g-2}(Y, \bbQ)^{-\epsilon} \quad \textnormal{and} \quad
 V_{+\epsilon} \;=\; H^{g-2}(X, \bbQ)^\perp \;\subseteq \; H^{g-2}(Y, \bbQ)^{+ \epsilon}
\]
where now $\perp$ denotes the orthocomplement inside the $+\epsilon$ eigenspace. Clearly~$V_-$ and $V_+$ are Hodge substructures of $H^{g-2}(Y, \bbQ)$. To compute their Hodge numbers one can use the Hirzebruch-Riemann-Roch theorem together with the following result, where $[\Theta]\in H^2(X, \bbQ)$ denotes the class of the theta divisor.

\begin{lem} \label{lem:chern-classes}
In terms of the restriction $\theta = [\Theta]|_Y\in H^2(Y, \bbQ)$, the Chern classes of~$Y$ are given by
\[ 
 c_i(Y) \;=\; (-1)^i \cdot (i+1) \cdot \theta^i \;\in \; H^{2i}(Y, \bbQ).
\] 
In particular, the variety $Y$ is of general type with canonical class $K_Y = 2\theta$.
\end{lem}

{\em Proof.} By definition we have $Y = \Theta \cap \Theta_{x}$ for some point $x\in X(\bbC)$. The embedding $Y\hookrightarrow \Theta$ gives an adjunction sequence
\[
0 \longrightarrow \calT_{Y} 
\longrightarrow \calT_{\Theta} \vert_{Y}
\longrightarrow \calO_X(\Theta_{x}) \vert_{Y}
\longrightarrow 0
\]
where $\calT_Y$ and $\calT_\Theta$ are the tangent bundles to $Y$ resp.~$\Theta$ (this makes sense also if the theta divisor has isolated singularities, provided that the point $x$ is chosen such that $Y$ does not meet the singularities). By restriction of the adjunction sequence for the embedding $\Theta \hookrightarrow X$ we also get an exact sequence
\[
 0\longrightarrow \calT_{\Theta}|_Y 
\longrightarrow \calT_X\vert_{Y}
\longrightarrow \calO_X(\Theta) \vert_{Y}
\longrightarrow 0.
\]
Since $\calT_X$ is trivial, it follows from these two adjunction sequences that the Chern polynomial $c_t(Y)= 1 + c_1(Y) \cdot t + c_2(Y) \cdot t^2 + \cdots$ satisfies 
$1 =c_t(Y) \cdot (1+\theta t)^2 $
because the Chern polynomial is multiplicative in short exact sequences. Now our claim about the Chern classes follows by a comparison of coefficients. The canonical sheaf $\omega_Y$ is by definition the top exterior power of the cotangent bundle $\calT_Y^*$, so for the canonical class we get
\[
 K_Y \;=\; c_1(\omega_Y) \;=\; c_1(\Lambda^{g-2}(\calT_Y^*)) \;=\; c_1(\calT_Y^*) \;=\; -c_1(Y) \;=\; 2\theta
\]
where in the third equality we have used that by the splitting principle the first Chern class of a bundle coincides with the first Chern class of any of its exterior powers~\cite[rem.~3.2.3(c)]{Fulton}. In particular, the canonical class $K_Y$ is ample and therefore also big in the sense of~\cite[sect.~2.2]{LazarsfeldPositivityI}. Thus $Y$ is of general type. \qed

\medskip

With notations as above, to compute intersection numbers between powers of the Chern classes of $Y$ we only need to know the image of $\theta^{g-2}$ under the degree isomorphism $\deg_Y: H^{g-2}(Y, \bbQ) \longrightarrow \bbQ$.
For this recall that $Y$ is the intersection of two general translates of the theta divisor. The intersection of cycles corresponds to the cup product on cohomology, so the fundamental class of~$Y$ is the cup product square $[Y] = [\Theta]^2 \in H^4(X, \bbQ)$. Hence
\[
 \deg_Y(\theta^{g-2}) \;=\; \deg_X([\Theta]^{g-2} \cdot [Y]) \;=\; \deg_X ( [\Theta]^g ) \;=\; g!
\]
where the last equality holds by the Poincar\'e formula~\cite[sect.~11.2.1]{BL}.

\begin{cor} \label{cor:gauss-bonnet}
For $g\geq 2$ the topological Euler characteristic of $Y$ is given by
\[\chi(Y) \;=\; (-1)^g\cdot (g-1)  \cdot g! \, . \] 
\end{cor}

{\em Proof.} By the Gauss-Bonnet theorem~\cite[p.~416]{GH} we have $\chi(Y)=\deg_Y(c_{g-2}(Y))$ so that the claim follows from lemma~\ref{lem:chern-classes} and from the above degree formula. \qed

\medskip

\section{Some low-dimensional examples}

To get a feeling for the above constructions, let us take a more explicit look at the situation for some small values of the dimension $g=\dim(X)$. We keep the notations from the previous section, always assuming that the theta divisor $\Theta \subset X$ is smooth or has only isolated singularities.

\medskip
 
{\em The case $g=2$}. Here $Y$ consists of two points $p,q$ which are interchanged by~ $\sigma$ as indicated in the following picture. 

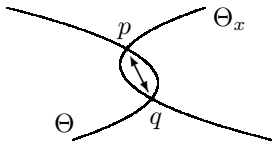
\begin{center}
\begin{picture}(150,65)

\put(70.5,39.25){\circle*{2}}
\put(80,20.4){\circle*{2}}

\put(67, 45){$p$}
\put(79, 11){$q$}

\put(72.5, 35){\vector(1,-2){6}}
\put(72.5, 35){\vector(-1,2){0.9}}

\put(43,8){$\Theta$}
\put(103,48){$\Theta_x$}

\qbezier(25,55)(125,30)(50,5)
\qbezier(125,5)(25,30)(100,55)

\end{picture}
\end{center}
 
\noindent
The fundamental classes $\alpha = [p]$ and $\beta = [q]$ of the two interchanged points form a basis of~$H^0(Y, \bbQ)$, and the eigenspace
\[ H^0(Y, \bbQ)^+ \;=\; \bbQ \cdot ( \alpha + \beta ) \;=\;  H^0(X, \bbQ) \]
equals the image of the weak Lefschetz embedding. Thus $V_+ = 0$, and the Hodge structure $V_- = \bbQ \cdot (\alpha - \beta)$ of weight zero has rank one. 
 
\medskip
 
{\em The case $g=3$}. Here $Y\rightarrow Y^+$ is an \'etale double covering of smooth projective curves, so we can consider the corresponding {\em Prym variety}. Before we come to our specific example, let us briefly recall the definition and some basic facts about Prym varieties in general, referring to~\cite{MuPrym} and~\cite[ch.~12]{BL} for details. For any finite covering $\tilde{C} \rightarrow C$ of smooth projective curves, the pushforward of divisors induces a norm epimorphism
$
 N_{\tilde{C}/C}:  J\tilde{C} \rightarrow JC
$
between the Jacobian varieties. For an \'etale double covering the kernel of this epimorphism has precisely two connected components, and one defines the Prym variety
\[
 P \;=\; \mathit{Prym}(\tilde{C}/C) \;=\; (\ker(N_{\tilde{C}/C}))^0 \;\subset\; J\tilde{C}
\]
to be the connected component which contains the origin. By construction this is an abelian subvariety of $J\tilde{C}$, and it turns out that the principal polarization on~$J\tilde{C}$ restricts to twice a principal polarization on~$P$. Prym varieties are useful in the study of ppav's since they are more general than Jacobian varieties but still accessible in terms of algebraic curves. Note that if the curve $C$ has genus $g+1$, then $\tilde{C}$ has genus $2g+1$ due to the Riemann-Hurwitz formula~\cite[cor.~IV.2.4]{Hartshorne}, and then the Prym variety has dimension 
\[ \dim(P) \;=\; (2g+1)-(g+1) \;=\; g. \]
Hence if we denote by $\calM_{g+1}$ the coarse moduli space of smooth projective curves of genus $g+1$, then for the coarse moduli space $\calR_{g+1}$ of \'etale double covers of such curves we have a diagram
\[
\xymatrix@R=0.8em{
 & \calR_{g+1} \ar[dl]_-\varphi \ar[dr]^-\pi \\
 \calM_{g+1} && \calA_{g}
}
\]
where $\varphi$ is the forgetful morphism and where $\pi$ is the morphism that assigns to a double covering its Prym variety. Here $\varphi$ is generically finite of degree $2^{g+1}-1$ since the \'etale double covers of a curve correspond bijectively to the \mbox{non-trivial} \mbox{two-torsion} points on its Jacobian variety. The Prym morphism $\pi$ has been studied a lot in relation with the moduli of abelian varieties. For $g\geq 6$ it is generically finite of degree one~\cite{FriedmanSmith}, though a construction of Donagi shows that it is never injective~\cite{Do}. In dimensions $g<6$ the Prym morphism $\pi$ is no longer generically finite but its fibres have a surprisingly rich geometric structure.

\medskip

This brings us back to the \'etale double covering $Y\rightarrow Y^+$ that is defined by the intersection of two general translates of the theta divisor on a complex ppav $X$ of dimension $g=3$. The curve $Y$ has genus $7$ by corollary~\ref{cor:gauss-bonnet}, so the associated Prym variety $P=\mathit{Prym}(Y/Y^+)$ is a complex ppav of dimension~$3$ as well. The following result goes back to the work of S.~Recillas~\cite{Re}, see also~\cite[\S 6.4]{Do}.

\begin{thm}
Let $X$ be a general complex ppav of dimension $g=3$.
\smallskip
\begin{enumerate}
 \item[\em (a)] For general $x\in X(\bbC)$ the Prym variety $P$ of the double cover $Y_x \rightarrow Y_x^+$ is isomorphic to the ppav $X$.
\smallskip
 \item[\em (b)] Every \'etale double covering of smooth curves with Prym variety isomorphic to $X$ arises as above for some point $x\in X(\bbC)$, and the coverings for two points $x_1, x_2$ are isomorphic iff $x_1 = \pm \, x_2$.
\end{enumerate}
\end{thm}

So via the above construction the fibre of the Prym morphism $\pi: \calR_4 \rightarrow \calA_3$ over a general ppav $X\in \calA_3(\bbC)$ is identified with a Zariski-open dense subset of the Kummer variety 
\[ K_X \;=\; X/\langle \pm \id_X\rangle. \] 
Fixing a general point $x\in X(\bbC)$, we now again write $Y=Y_x$ etc. The definition of the Prym variety $P$ as a component of the kernel of the norm epimorphism shows that the Jacobian variety $JY$ is isogenous to $P\times JY^+$, which together with part~{\em (a)} of the theorem implies
\[ H^1(Y, \bbQ) \; = \; H^1(X,\bbQ)\oplus H^1(Y^+,\bbQ). \] 
Hence in this case $V_- = 0$, and the eigenspace $V_+ = H^1(Y^+, \bbQ)$ is a pure Hodge structure of weight one and rank eight.
 
\medskip
 
\label{subsec:g4}
{\em The case $g=4$}. Here $Y\rightarrow Y^+$ is an \'etale double covering of smooth projective surfaces that will occupy us for the rest of this paper. 
We claim that $V_-$ has Hodge numbers $h^{2,0}=h^{0,2}=11$ and $h^{1,1}=30$ whereas $V_+$ is of pure Hodge type $(1,1)$ and rank six. Indeed, bearing in mind that
\[ H^2(Y^+, \bbQ) \;=\; H^2(Y, \bbQ)^+ \;=\; H^2(X, \bbQ) \oplus V_+ \] 
this is a direct consequence of the numerical data in the following table together with the Hodge numbers $h^{2,0}(X)=h^{0,2}(X)=6$ and $h^{1,1}(X)=16$.

\begin{lem} \label{lem:hodge-numbers}
For a complex ppav $X$ of dimension $g=4$ whose theta divisor has at most isolated singularities, we have the following Hodge numbers.
\smallskip
\[
\begin{array}{l|clc}
  & \;\; h^{2,0} = h^{0,2} \;\; & \;\; h^{1,1} \;\; & \;\; h^{1,0} = h^{0,1} \;\; \\ \hline 
 Y & 17 & 52 & 4\\ 
 Y^+ & 6 & 22 & 0 
\end{array}
\medskip
\]
\end{lem}

{\em Proof.} The last column is clear since $H^1(Y, \bbQ)=H^1(Y, \bbQ)^-=H^1(X, \bbQ)$ by the weak Lefschetz theorem. Furthermore from corollary~\ref{cor:gauss-bonnet} we obtain the topological Euler characteristic $\chi(Y)= 72$, which implies~$h^2(Y) = 86$ since $h^0(Y)=h^{4}(Y)=1$ and $h^1(Y)=h^3(Y)=8$. Lemma~\ref{lem:chern-classes} together with the Hirzebruch-Riemann-Roch theorem also gives the holomorphic Euler characteristic
\[ 
\chi(\calO_{Y})
\;=\;
\deg_Y(c_1^2(Y) + c_2(Y)) / 12
\;=\;
\deg_Y((-2\theta)^2 + 3\theta^2) / 12
\;=\;
14.
\]
Plugging in the values $h^{0,0}(Y)=1$ and $h^{0,1}(Y)=4$ we get $h^{0,2}(Y)=17$, and this implies $h^{1,1}(Y)=52$ because the Hodge numbers in degree two must sum up to the second Betti number $h^2(Y)=86$.

\medskip

To obtain also the Hodge numbers of the quotient surface $Y^+=Y/\sigma$ we use that the quotient morphism $q: Y\to Y^+$ is an \'etale double covering. In particular we have $\deg_Y(q^*(\alpha)) = 2 \deg_{Y^+}(\alpha)$ for all $\alpha \in H^4(Y^+, \bbQ)$, and $c_i(Y) = q^*(c_i(Y^+))$ for all~$i\in \bbN$. So the Gauss-Bonnet and Hirzebruch-Riemann-Roch theorems show 
\[ \chi(Y^+) \;=\;\chi(Y)/2 \;=\; 36 \quad \textnormal{and} \quad \chi(\calO_{Y^+})\;=\;\chi(\calO_Y)/2\;=\;7. \]
Hence it follows that $h^2(Y^+)=36-1-1=34$ and $h^{0,2}(Y^+)=7-1=6$ by similar computations as above, and we are done. \qed

\medskip

As a corollary we obtain that the intersection form has the following signatures on the various subspaces of $H^2(Y, \bbQ)$, where we denote by $s_+$ and $s_-$ the dimension of maximal subspaces on which the form is positive resp.~negative definite.

\begin{cor} \label{cor:signatures}
For $g=4$ the intersection form has the following signatures.
\[
\begin{array}{r|ccccc}
  & H^2(Y, \bbQ) & H^2(Y^+, \bbQ) & H^2(X, \bbQ) & \;\;\;\; V_- \;\;\;\; & \;\;\;\; V_+ \;\;\;\; \\ \hline 
  s_+ & 35  & 13 & 13 & 22 & 0 \\ 
   s_- & 51  & 21 & 15 & 30 & 6 
\end{array}
\]
\end{cor}

{\em Proof.} By the Hodge index theorem~\cite[th.~6.33]{VoisinHodge1}, for a smooth compact K\"ahler surface $S$ with K\"ahler class $\omega$ the intersection form is 
\begin{itemize}
\item positive definite on $(H^{2,0}(S)\oplus H^{0,2}(S))\cap H^2(S, \bbR)$ and on $\bbR\cdot \omega$,
\item negative definite on the orthocomplement of  $\omega$ in $H^{1,1}(S)\cap H^2(S, \bbR)$.
\end{itemize}
Hence our claim follows by taking $S=Y$ and using lemma~\ref{lem:hodge-numbers}. \qed

\medskip

In particular, on the subspace $V_+$ the intersection form is negative definite. In the next two sections we will work on the level of integral cohomology and determine the lattices whose signatures have been listed above. 

\medskip

\section{Integral cohomology} \label{sec:coh}

Let $X$ be a complex ppav of dimension $g=4$ and $\Theta \subset X$ a smooth symmetric theta divisor. Fix a general point $x\in X(\bbC)$ and consider the \'etale double covering of smooth surfaces 
\[
  Y \;=\; \Theta \cap \Theta_x \;\longrightarrow \;  Y^+ = Y/\sigma
\]
for the covering involution $\sigma=\sigma_x$ as defined in corollary~\ref{cor:smooth-and-etale}. We claim that the integral cohomology of our surfaces looks as follows.

\begin{prop} \label{prop:coh}
The smooth projective surfaces $Y$ and $Y^+$ have the following integral cohomology groups.
\[
\begin{array}{r|ccccc} 
 n \; &  \quad 0 & \qquad 1 \; & 2 & 3 & \quad 4  \\ \hline 
 H^n(Y, \bbZ) \; & \quad \bbZ & \qquad \bbZ^8 \; & \bbZ^{86} & \bbZ^8 & \quad \bbZ \rule{0em}{1.1em} \\ 
 H^n(Y^+, \bbZ) \; & \quad \bbZ & \qquad 0 \; & \bbZ^{34} \oplus (\bbZ/2\bbZ)^9 & (\bbZ/2\bbZ)^9 & \quad \bbZ
\end{array}
\]
Moreover, the pull-back map for the quotient morphism $q:Y \longrightarrow {Y^+}$ gives rise to an epimorphism onto the $\sigma$-invariants
\[ q^*:  \; \; 
 H^2({Y^+}, \bbZ) \; \twoheadrightarrow \; H^2(Y, \bbZ)^+ \; \subset \; H^2(Y, \bbZ) 
\]
and the kernel of this epimorphism is equal to the torsion subgroup of $H^2({Y^+}, \bbZ)$. 
\end{prop}

{\em Proof.} 
For any complex projective variety $V$ of dimension $d+1$ with an ample effective divisor $W\subset V$ such that the complement $V\setminus W$ is smooth, the weak Lefschetz theorem says that the restriction and pushforward maps
\[
 H^n(V, \bbZ) \; \longrightarrow \; H^n(W, \bbZ)
 \quad \textnormal{and} \quad
 H_n(W, \bbZ) \; \longrightarrow \;  H_n(V, \bbZ)
\]
are isomorphisms in degree~$n<d$ and a mono- resp.~epimorphism in degree~$n=d$, see~\cite[cor.~7.3]{MilnorMorse} and~\cite[th.~1.23]{VoisinHodge2}.
Applying this first to the divisor $\Theta \subset X$ and then to the divisor $Y \subset \Theta$ we obtain that
\[
 H^n(Y, \bbZ) \;\cong\; H^n(X, \bbZ) \quad \textnormal{and} \quad
 H_n(Y, \bbZ) \;\cong\; H_n(X, \bbZ) \quad \textnormal{for} \quad 
 n<2. 
\]
In particular, the weak Lefschetz theorem for cohomology implies that $H^n(Y, \bbZ)$ has the claimed form for $n<2$. For $n>2$ the corresponding statement follows from the weak Lefschetz theorem for homology because $H^n(Y, \bbZ)\cong H_{4-n}(Y, \bbZ)$ by Poincar\'e duality. For any~$n$ we furthermore have an exact sequence 
\vspace{0.5em}
\[
 0 
 \; \longrightarrow \; \underbrace{\mathit{Ext}(H_{n-1}(Y, \bbZ), \bbZ)}_{\textnormal{torsion group}} 
 \; \longrightarrow \; H^n(Y, \bbZ) 
 \; \longrightarrow \; \underbrace{\Hom(H_n(Y, \bbZ), \bbZ)}_{\textnormal{torsion-free}} 
 \; \longrightarrow \; 0 
 \vspace{-0.1em}
\]
by the universal coefficient theorem. Using that $H_1(Y, \bbZ) \cong H_1(X, \bbZ)$ is torsion-free so that $\mathit{Ext}(H_1(Y, \bbZ), \bbZ)=0$, we obtain that
\[
 H^2(Y, \bbZ) \;\cong \; \Hom(H_2(Y, \bbZ), \bbZ)
\]
is a free abelian group, and the rank of this group is $h^2(Y) = 86$ by lemma~\ref{lem:hodge-numbers}.

\medskip

It remains to show that $q^*: H^2(Y^+, \bbZ) \longrightarrow H^2(Y, \bbZ)$ is an epimorphism onto the $\sigma$-invariants whose kernel is equal to the torsion subgroup, and to determine the groups $H^n(Y^+, \bbZ)$. Note that
$
 H_n(Y^+, \bbZ) \cong H^{4-n}(Y^+, \bbZ) 
$
for all $n$ by Poincar\'e duality, so by the universal coefficient theorem we have (non-canonically)
\begin{align} \nonumber
 H^n(Y^+, \bbZ) 
 & \;\; \cong \;\; 
 H_{n-1}(Y^+, \bbZ)_{\mathit{tor}}  \oplus  H_n(Y^+, \bbZ)_{\mathit{free}} \\ \nonumber
 & \;\; \cong \;\;
 H^{5-n}(Y^+, \bbZ)_{\mathit{tor}}  \oplus  H^{4-n}(Y^+, \bbZ)_{\mathit{free}}
\end{align}
where the subscripts {\em tor} and {\em free} refer to the maximal torsion subgroup and the maximal free abelian subgroup, respectively. So for the computation of $H^n(Y^+, \bbZ)$ it suffices to deal with degrees $n\leq 2$. Now we can argue as follows.

\medskip

Recall from topology~\cite[ex.~VII.4(b)]{BrownCohomology} that if $W\longrightarrow W^+=W/\Gamma$ is any regular covering map between connected, locally pathwise connected topological spaces with Galois group $\Gamma = Aut(W/W^+)$, then one has a spectral sequence 
\[ E_2^{pq} \;=\; H^p(\Gamma, H^q(W, \bbZ)) \; \Longrightarrow \; H^{p+q}(W^+, \bbZ) \]
where the terms on the left hand side denote group cohomology with coefficients in the $\Gamma$-module $H^q(W, \bbZ)$. We apply this for $W=Y$ and $\Gamma = \langle \sigma \rangle$. To compute some of the $E_2$-terms in this case, recall that the group cohomology of any module~$M$ under the cyclic group $\langle \sigma \rangle \cong \bbZ/2\bbZ$ is 
\[
 H^p(\langle \sigma \rangle, M) \;=\; 
\begin{cases}
  \mathrm{ker}(\sigma - 1) & \textnormal{for $p=0$}, \\
  \mathrm{ker}(\sigma + 1) / \mathrm{im}(\sigma - 1) & \textnormal{for $p>0$ odd}, \\
  \mathrm{ker}(\sigma - 1) / \mathrm{im}(\sigma + 1) & \textnormal{for $p>0$ even}. 
\end{cases}
\]
Let us plug in $M=H^q(Y, \bbZ)$ with $q\in \{0, 1, 2\}$. Clearly on $H^0(Y, \bbZ)\cong \bbZ$ the involution $\sigma$ acts trivially, whereas by the weak Lefschetz theorem we know that on $H^1(Y, \bbZ) \cong H^1(X, \bbZ) \cong \bbZ^8$ it acts by $-1$. Hence 
\begin{align} \nonumber
 \frac{\mathrm{ker}(\sigma + 1 \mid H^q(Y, \bbZ))}{\mathrm{im}(\sigma - 1\mid H^q(Y, \bbZ))} 
 \;\; & \cong \;\; 
 \begin{cases}
 0 & \textnormal{for} \; \; q = 0, \\
 (\bbZ/2\bbZ)^8 & \textnormal{for} \; \; q = 1,
 \end{cases}
\\ \nonumber
 \frac{\mathrm{ker}(\sigma - 1 \mid H^q(Y, \bbZ))}{\mathrm{im}(\sigma + 1\mid H^q(Y, \bbZ))} 
 \;\; &\cong \;\; 
 \begin{cases}
 \bbZ/2\bbZ \;\;\; \;\; \, & \textnormal{for} \; \; q = 0, \\
 0 & \textnormal{for} \; \; q = 1.
 \end{cases}
\end{align}
In other words, in the case at hand the $E_2$-tableau of the spectral sequence takes the following form.
\[
\xymatrix@C=1em@R=0.7em{
q & \vdots & \vdots & \vdots & \vdots &  \\
2 & H^2(Y, \bbZ)^+ \ar[drr] & \ast \ar[drr] & \ast & \ast & \cdots &  \\
1 & 0 \ar[drr] & (\bbZ/2\bbZ)^8  \ar[drr] & 0 & (\bbZ/2\bbZ)^8 & \cdots & \\
0 & \bbZ & 0 & \bbZ/2\bbZ & 0 & \cdots & \\
 \ar@<2ex>[rrrrrr] \ar@<-2ex>[uuuu] & 0 & 1 & 2 & 3 & & p
}
\]
Hence $H^0(Y^+, \bbZ)=\bbZ$ and $H^1(Y^+, \bbZ)=0$. The zeroes for $(p,q) = (2,1), (3,0)$ also show that for the graded object with respect to the limit filtration we have
\[
 Gr_i(H^2(Y^+, \bbZ)) \;=\; 
 \begin{cases}
 H^2(Y, \bbZ)^+ & \textnormal{for $i=2$}, \\
 \hspace{0.25em} (\bbZ/2\bbZ)^8 & \textnormal{for $i=1$}, \\
 \hspace{0.6em} \bbZ/2\bbZ & \textnormal{for $i = 0$},
 \end{cases}
\]
and the top quotient defines a surjection
\[
 H^2({Y^+}, \bbZ) \;=\; E^2_\infty \; \twoheadrightarrow \; E_\infty^{02} \;=\; E_2^{02} \;=\; H^2(Y, \bbZ)^+.
\]
By construction of the spectral sequence this surjection coincides with the map $q^*$ induced by the quotient morphism $q$. The spectral sequence also shows that the kernel $K$ of this surjection is an extension of $(\bbZ/2\bbZ)^8$ by $\bbZ/2\bbZ$. In fact $K=(\bbZ/2\bbZ)^9$ because the kernel $K$ must be a $2$-torsion group: If we denote by 
\[ q_!: \quad H^2(Y, \bbZ) \; \longrightarrow \; H^2(Y^+, \bbZ) \] 
the Gysin map (the Poincar\'e dual to the pushforward on homology), then $q_! q^*$ is multiplication by the degree $\deg(q)=2$.
To conclude the proof, we observe that the group $F=H^2(Y^+, \bbZ)/K$ is a subgroup of the free abelian group $H^2(Y, \bbZ)$. As such it is torsion-free, and its rank is $h^2(Y^+)=34$ by lemma~\ref{lem:hodge-numbers}.
\qed

\medskip

\section{N\'eron-Severi lattices} \label{sec:lattice}

Recall that in the setting of the previous section we have embeddings of free abelian groups
\[
 H^2(X, \bbZ) \; \subset \; H^2(Y, \bbZ)^+ \; \subset \; H^2(Y, \bbZ).
\]
In what follows we equip these groups with the integral bilinear form that is given by the intersection form on $Y$. The goal of the present section is to determine the structure of the obtained integral lattices with respect to this bilinear form. Let us briefly recall some basic notions from lattice theory~\cite{CS2}.

\medskip

By a {\em lattice} we mean a pair consisting of a finitely generated free abelian group~$L$ and a non-degenerate symmetric bilinear form $b: L\times L \longrightarrow \bbZ$, though we will usually suppress the bilinear form in the notation. Such a lattice is called even if~$b(\lambda, \lambda)\in 2\bbZ$ for all $\lambda\in L$. For integers $n\neq 0$ we denote by $L(n)$ the lattice with the same underlying free abelian group as~$L$ but with the bilinear form $n\cdot b$ in place of $b$. As usual, by an embedding of lattices we mean an embedding of free abelian groups which respects the bilinear forms, and similarly for isomorphisms of lattices. For example, the diagonal embedding 
\[ \diag: \quad L(2) \; \hookrightarrow \; L^2 \; = \; L \oplus L \]
is given by $\lambda\mapsto (\lambda, \lambda)$. We also remark that any lattice can be defined with respect to a chosen basis $e_1, \dots, e_n$ of the underlying free abelian group by the associated {\em Gram matrix} --- the non-degenerate integral symmetric matrix $(b(e_i, e_k))_{i,k = 1, \dots, n}$ of size $n\times n$ whose entries are the scalar products between the basis vectors. The {\em discriminant} of a lattice is defined as the determinant of its Gram matrix with respect to any chosen basis, and the lattice is called {\em unimodular} if its discriminant is $\pm 1$. In what follows we consider the lattice
\[
\Lambda \;=\; \bbZ^4 \quad \textnormal{with Gram matrix} \quad 
\left(
\begin{smallmatrix}
0&1&1&1\\ 1&0&1&1\\ 1&1&0&1\\ 1&1&1&0
\end{smallmatrix}
\right)
\]
and the standard hyperbolic lattice $U=\bbZ^2$ with Gram matrix $\left( \begin{smallmatrix} 0 & 1\\ 1 & 0 \end{smallmatrix} \right)$. Both $\Lambda$ and~$U$ are even lattices. Note that~$U$ is unimodular whereas a direct calculation shows that $\Lambda$ has discriminant $-3$. This being said, in terms of the even lattices
\[
 K \;=\; U^{12}\oplus \Lambda \quad \textnormal{and} \quad L \;=\; U^{13}\oplus E_8(-1).
\]
the main result of the present section can be formulated as follows.

\begin{prop}  \label{prop:E6} \label{prop:summary}
There exists an isomorphism $H^2(Y, \bbZ) \cong L^2 \oplus U^9$ of lattices
which gives a commutative diagram
\[
\xymatrix@C=3.5em@M=0.8em{
H^2(X, \bbZ) \ar@{^{(}->}[r] \ar[d]^-\cong
& H^2(Y, \bbZ)^+ \ar@{^{(}->}[r] \ar[d]^-\cong 
& H^2(Y, \bbZ) \ar[d]^-\cong \\
K(2)  \ar@{^{(}->}[r]^-{i(2)}
& L(2)   \ar@{^{(}->}[r]^-{(\diag, 0)} 
& L^2 \oplus U^9
}
\]
for some lattice embedding $i: K\hookrightarrow L$ with orthocomplement $K^\perp \cong E_6(-1)$. 
\end{prop}

To complete the picture, we remark that with the above identifications it follows that on the lattice $L^2\oplus U^9$ the involution $\sigma: Y\longrightarrow Y$ acts by the formula 
\[ \sigma(\lambda_1, \lambda_2, u) \;=\; (\lambda_2, \lambda_1, -u) 
\quad \textnormal{for} \quad u\in U^9 \quad \textnormal{and} \quad  \lambda_1, \lambda_2\in L.
\]
Indeed, looking at the eigenspace decomposition over the rational numbers one sees that $\sigma$ must act by $-1$ on the orthocomplement of $H^2(Y, \bbZ)^+$ in~$H^2(Y, \bbZ)$.

\medskip

The proof of proposition~\ref{prop:E6} will occupy the rest of this section and will also involve a study of the intersection form on $H^2(Y^+, \bbZ)$. Note that since this latter group contains torsion, it is not a lattice in the above sense. So let us introduce the following notations. For any smooth complex projective surface $S$ let
\[
 L_S \;=\; H^2(S, \bbZ) / \mathit{Torsion}
\]
be the quotient of the finitely generated abelian group $H^2(S, \bbZ)$ by its torsion subgroup. Since the intersection form
\[
 b: \quad H^2(S, \bbZ) \times H^2(S, \bbZ) \; \longrightarrow \; \bbZ
\]
takes its values in the torsion-free group $\bbZ$, the intersection pairing between any torsion class and any other class vanishes. Hence the intersection form of the surface factors over a non-degenerate symmetric bilinear form $b: L_S \times L_S \longrightarrow \bbZ$, and we view $L_S$ as a lattice with respect to this bilinear form. By Poincar\'e duality this lattice is always unimodular. We also have the following well-known

\begin{lem} \label{lem:even-lattice}
If for some integral class $\alpha \in H^2(S, \bbZ)$ the first Chern class satisfies the congruence $c_1(S) \equiv 2\alpha$ modulo torsion, then $L_S$ is an even lattice. 
\end{lem}

{\em Proof.} By naturality the reduction homomorphism $r: \bbZ \longrightarrow \bbZ/2\bbZ$ gives rise to a commutative diagram
\[
 \xymatrix@C=3em@M=0.6em@R=1.8em{
 H^2(S, \bbZ) \times H^2(S, \bbZ) \ar[r]^-\cup \ar[d]_-{(r_*, r_*)} & H^{4}(S, \bbZ) \ar[r]^-{\deg_S} \ar[d]_-{r_*} & \bbZ \ar[d]^-r \\
 H^2(S, \bbZ/2) \times H^2(S, \bbZ/2) \ar[r]^-\cup & H^{4}(S, \bbZ/2) \ar[r]^-{\deg_S}  & \bbZ/2
 }
\]
%
The composite of horizontal arrows in the upper row is the intersection pairing we are interested in. Now the image of the first Chern class $c_1(S)$ under the reduction map
$r_*: H^2(S, \bbZ) \longrightarrow  H^2(S, \bbZ/2\bbZ)$
is the Stiefel-Whitney class $w_2(S)$ \cite[p.~171]{MS}, and for the intersection form modulo two we have
\[
 \deg_S(w_2(S) \cup r_*(\beta)) \; = \; r(\deg_S(\beta^2)) \quad \textnormal{for all} \quad \beta \in H^2(S, \bbZ) 
\]
by the Wu formula~\cite[prop.~1.4.18]{GompfStipsicz}. Hence the claim follows. \qed

\medskip

Let us now come back to the smooth surfaces $Y$ and $Y^+$ which are defined by a general intersecion of translates of a smooth theta divisor on a complex ppav $X$ of dimension $g=4$. We have the following abstract isomorphisms.

\begin{lem} \label{rem:Ylattice}
With notations as above, 
\[ L_Y \;\cong \;  L^2 \oplus U^9 \quad \textnormal{\em and} \quad \; L_{Y^+} \cong  L. \]
\end{lem}

{\em Proof.} By corollary~\ref{cor:signatures} the intersection form on $L_Y$ and $L_{Y^+}$ is indefinite of signature $(35, 51)$ resp.~$(13, 21)$. An even unimodular lattice is determined uniquely by its rank and signature~\cite{SeFormesBilineaires}, so it only remains to show that both $L_Y$ and $L_{Y^+}$ are even lattices. For this we can apply lemma~\ref{lem:even-lattice}. Indeed, the condition of the lemma holds for $L_Y$ because 
$c_1(Y) = 2\alpha$ for the class $\alpha = -[\Theta]|_Y \in H^2(Y, \bbZ)$
by lemma~\ref{lem:chern-classes}. Furthermore, by proposition~\ref{prop:coh} we can find a vector $\alpha^+ \in L_{Y^+}$ such that $q^*(\alpha^+)=\alpha$ where $q: Y\rightarrow Y^+$ is the quotient morphism. Then
\[
 q^*(c_1(Y^+)-2\alpha^+ )  \;=\; c_1(Y) - 2\alpha \;=\; 0
\]
so that again by proposition~\ref{prop:coh} the class $c_1(Y^+)-2\alpha^+$ is a torsion class. Hence lemma~\ref{lem:even-lattice} shows that $L_{Y^+}$ is an even lattice as well. \qed

\medskip

Now that we have found the abstract isomorphism type of the lattices~$L_Y$ and~$L_{Y^+}$, the next thing to be done is to determine the relationship between these two lattices. For this we denote by
\[
 L_Y^+ \; = \; H^2(Y, \bbZ)^+ \; \subset \; L_Y \;=\; H^2(Y, \bbZ)
\]
the invariants of $L_Y$ under the action of the involution $\sigma$, i.e.~the lattice which occurs in the middle part of the diagram in proposition~\ref{prop:E6}.

\begin{lem} \label{lem:quotient} 
The pull-back under the quotient morphism $q: Y\longrightarrow Y^+$ gives rise to a lattice isomorphism
\[ q^*: \;\; L_{Y^+}(2) \; \stackrel{\cong}{\longrightarrow} \; L_Y^+ \; \subset \; L_Y \,.
\medskip \]
\end{lem}

{\em Proof.} By proposition~\ref{prop:coh} the pull-back $q^*$ gives an isomorphism of the underlying abelian groups. Via this isomorphism, the intersection form on $L_Y$ induces twice the intersection form on $L_{{Y^+}}$ since $q: Y\rightarrow {Y^+}$ is an \'etale double cover: The degree maps for the top cohomology of our surfaces satisfy 
\[ \deg_Y( q^*(a)\cup q^*(b) ) \;=\; \deg _Y(q^*(a \cup b)) \;=\; 2 \deg_{Y^+} (a\cup b) \]
for all cohomology classes~$a, b$ in~$H^2({Y^+}, \bbZ)$. \qed

\medskip

The above two lemmas give the middle and the right vertical isomorphism in the diagram of proposition~\ref{prop:E6}. These isomorphisms can be chosen such that the square on the right hand side of the diagram commutes, indeed we have the following

\begin{lem}
Let $L_1, L_2 \subset L_Y$ be sublattices of rank $r\leq 34$ which are primitive in the sense that the quotients $L_Y/L_1$ and $L_Y/L_2$ are torsion-free. Then any lattice isomorphism $\varphi_0: L_1 \rightarrow L_2$ extends to a lattice automorphism $\varphi: L_Y \rightarrow L_Y$.
\end{lem}

{\em Proof.} This follows from the criterion for the extension of lattice isomorphisms given in~\cite{Ja}. Let us briefly check the two conditions in loc.~cit. The first condition is that the rank $r$ of the two primitive sublattices and the signature $(s_+, s_-)=(35,51)$ of the ambient lattice $L_Y$ satisfy 
\[ 2(r + 1) \; \leq \; (s_+ + s_-)-|s_+-s_-|,
\]
which in our case holds for $r\leq 34$. The second condition concerns primitive vectors $v\in L_Y$ which are characteristic in the sense that  $b(v, w) \equiv b(w, w) \; \mathrm{mod} \; 2\bbZ$ for all vectors $w\in L_Y$. But in the case at hand this condition is void because an even unimodular lattice does not contain any primitive characteristic vectors.  \qed

\medskip

To finish the proof of proposition~\ref{prop:E6} it remains to deal with the square on the left hand side in the diagram of the proposition. For this we need to compute the image of the weak Lefschetz embedding, i.e.~the sublattice 
\[ L_X \;=\; H^2(X, \bbZ) \;\subset \; L_Y^+ \;=\; H^2(Y, \bbZ)^+ . \] 
Note that $K = U^{12} \oplus \Lambda$ has signature $(13, 15)$ and discriminant~$-3$ whereas the even unimodular lattice $L=U^{13}\oplus E_8(-1)$ has signature $(13, 21)$. Accordingly the orthocomplement $K^\perp$ of any embedding 
\[ i: \;\; K \; \hookrightarrow \; L \]
is a negative definite even lattice of rank six with discriminant $3$, hence isomorphic to $E_6(-1)$ by the classification of lattices with small discriminant in~\cite{CS}. Thus the proof of proposition~\ref{prop:E6} is completed by the following computation.

\begin{lem} \label{lem:Xlattice}
With notations as above, there exists an isomorphism of lattices 
\[ L_X \; \stackrel{\cong}{\longrightarrow} \; K(2). \]
\end{lem}

{\em Proof.} 
Consider the embedding $i: Y = \Theta \cap \Theta_x \hookrightarrow X$. By definition the Gysin morphism $i_!$
on cohomology is the Poincar\'e dual of the pushforward morphism $i_*$ on homology, so we have a commutative diagram
\[
\xymatrix@M=0.5em@C=3em@R=1.8em{
 H^n(Y, \bbZ) \ar[r]^-{i_!} \ar[d]_-\cong & H^{n + 4}(X, \bbZ) \ar[d]^-\cong \\
 H_{4-n}(Y, \bbZ) \ar[r]^-{i_*} & H_{4-n}(X, \bbZ).
}
\]
If~$\mathbf{1}\in H^0(Y, \bbZ)$ denotes the unit element of the cohomology ring of the surface $Y$, the fundamental cohomology class of~$Y$ in $X$ is $i_!(\mathbf{1}) = [Y] \in H^4(X, \bbZ)$. So the projection formula shows
\[
  i_! ( i^* (\gamma)) \;=\; i_! ( i^* (\gamma)\cup \mathbf{1}) \;=\; \gamma \cup i_!(\mathbf{1}) \;=\; \gamma \cup [Y] \quad 
\]
for all $\gamma \in H^4(X, \bbZ)$. Now it follows from the commutative diagram
\[
\xymatrix@C=3em@M=0.5em@R=1.8em{
 H^2(X, \bbZ) \times H^2(X, \bbZ) \ar[r]^-{\cup} \ar[d]^-{i^*} & H^4(X, \bbZ) \ar[r]^-{-\cup [Y]} \ar[d]^-{i^*} & H^8(X, \bbZ) \ar[d]^-{\deg_X} \\
 H^2(Y, \bbZ) \times H^2(Y, \bbZ) \ar[r]^-{\cup} & H^4(Y, \bbZ) \ar[r]^-{\deg_Y} \ar[ur]^-{i_!} & \bbZ \\
}
\]
that the intersection form on $H^2(Y, \bbZ)$ restricts on the subspace $H^2(X, \bbZ)$ to the bilinear form
\[
 b: \; H^2(X, \bbZ) \times H^2(X, \bbZ) \longrightarrow \bbZ, \; (\alpha, \beta) \mapsto \deg_X(\alpha \cup \beta \cup [Y]).
\]
We want to determine this bilinear form. The intersection product in homology corresponds to the cup product in cohomology, so the fundamental class of $Y$ in cohomology is
\[
  [Y] \;=\; [\Theta \cap \Theta_{x}] \;=\; [\Theta] \cup [\Theta_{x}] \;=\; [\Theta]^2.
\]
To compute this class explicitly, we choose a basis of the cohomology ring of $X$ as follows. The principal polarization on $X$ is given by an alternating bilinear form on~$H_1(X, \bbZ)$. Take an integral basis 
\[ \lambda_1, \dots, \lambda_4, \mu_1, \dots, \mu_4 \; \in \;  H_1(X, \bbZ) \; \cong \;  \bbZ^8 \]
in which this bilinear form is given by the scalar products $(\lambda_i, \lambda_k)=(\mu_i, \mu_k)=0$ and $(\lambda_i, \mu_k)=-(\mu_k, \lambda_i) = \delta_{ik}$ for all $i$, $k$. Let
\[ x_1, \dots, x_4, y_1, \dots, y_4
 \;\in \; H^1(X, \bbZ) \; = \; \Hom(H_1(X, \bbZ), \bbZ)
\]
be the dual basis. Since the cohomology class $[\Theta ]=c_1(\calO_X(\Theta))$ defines the principal polarization, it is given in the exterior algebra $H^\bullet(X, \bbZ)=\Lambda^\bullet(H^1(X, \bbZ))$ by
\[ [\Theta] \;=\; \sum_{i=1}^4 \, x_i\wedge y_i \; \in \; H^2(X, \bbZ). \]
The cup product in cohomology corresponds to the wedge product in the exterior algebra, so we get
\[
  [Y] \;=\; [\Theta]^2 \;=\; \Bigl(\sum_{i=1}^4 x_i\wedge y_i \Bigr) \wedge \Bigl(\sum_{j=1}^4 x_j\wedge y_j\Bigr) \;=\;
  2\cdot\hspace{-1em} \sum_{1\leq i < j \leq 4} x_i\wedge y_i \wedge x_j \wedge y_j.
\]
In particular, $v\wedge [Y] = 0$ for all vectors $v\in H^4(X, \bbZ)$ of the form $v = x_k \wedge x_l\wedge x_m \wedge x_n$ or $v=x_k\wedge x_l\wedge x_m\wedge y_n$ or $v=x_k\wedge y_l\wedge y_m\wedge y_n$ or~$v = y_k \wedge y_l\wedge y_m \wedge y_n$ with arbitrary indices $k, l, m, n$, and also for all vectors of the form $v=x_k\wedge x_l\wedge y_m\wedge y_n$ with $\{k, l\} \neq \{m, n\}$. However, for $k \neq l$ one easily checks that
\[ x_k\wedge x_l \wedge y_k \wedge y_l \wedge [Y] \;=\; -2\cdot \omega
\]
for the class
$\omega = x_1\wedge y_1 \wedge \cdots \wedge x_4\wedge y_4 = \tfrac{1}{4!} \cdot [\Theta]^4$ 
generating $H^8(X, \bbZ)$.

\medskip

Now consider the basis of the lattice $L_X = \Lambda^2(H^1(X, \bbZ))$ consisting of the vectors $u_{ik} = x_i\wedge y_k$ with $i, k\in \{ 1, 2, 3, 4\}$ and of the vectors $v_{ik}=x_i\wedge x_k$ and~$w_{ik}=y_i\wedge y_k$ with $1\leq i < k \leq 4$. By the above the bilinear form $b$ is non-zero only on those pairs of the above vectors whose wedge product has the form $\pm \hspace{0.1em} x_l \wedge x_m \wedge y_l \wedge y_m$ with $l \neq m$. So the only non-zero scalar products between our basis vectors are 
\smallskip
\begin{itemize}
\item[]
$b(u_{ii}, u_{kk})=x_i\wedge y_i\wedge x_k\wedge y_k \wedge [Y] = 2$ for $i\neq k$, 
\item[]
$b(u_{ik}, u_{ki})=x_i\wedge y_k \wedge x_k\wedge y_i\wedge [Y] = -2$ for $i\neq k$, and 
\item[]
$b(v_{ik}, w_{ik})=b(w_{ik}, v_{ik}) = y_i\wedge y_k \wedge x_i\wedge x_k \wedge [Y] = -2$ for $i<k$. 
\end{itemize}
\smallskip
Hence we obtain an orthogonal sum decomposition \medskip
\[
 L_X \;=\; 
 \langle u_{11}, \dots, u_{44} \rangle \; \; \oplus
 \bigoplus_{1\leq i < k \leq 4} \langle u_{ik}, u_{ki} \rangle 
 \oplus \langle v_{ik}, w_{ik} \rangle 
\]
where $\langle u_{11}, \dots, u_{44} \rangle \cong \Lambda(2)$ and where $ \langle u_{ik}, u_{ki} \rangle \cong  \langle v_{ik}, w_{ik} \rangle  \cong U(2)$. \qed

\medskip

\section{Variations of Hodge structures} \label{sec:translates}

So far we have always fixed a general point $x\in X(\bbC)$. Now let us see what happens when this point varies. For the moment the dimension $g=\dim(X)$ can be arbitrary, but as before we always assume that the theta divisor of our ppav is smooth or has only isolated singularities. To begin with, we put the rational Hodge structures $V_\pm$ from section~\ref{sec:intersections} into a family as follows.

\begin{lem} \label{lem:vhs}
Over some Zariski-open dense $U\subset X$ there exists for each $n\geq 0$ a polarized variation $\bbH^n$ of pure Hodge structures  with stalks
\[ \bbH^n_x \;=\; H^n(Y_x, \bbQ), \]
and subvariations
$\bbV_\pm \subset \bbH^{g-2}$ whose stalks are the subspaces $V_\pm \subset H^{g-2}(Y_x, \bbQ)$.
\end{lem}

{\em Proof.} Consider the Zariski-open dense subset $U\subset X$ and the smooth proper family $f_U: Y_U \rightarrow U$ in lemma~\ref{lem:bertini} and corollary~\ref{cor:smooth-and-etale}. We can put
$
  \bbH^n = R^nf_{U*}(\bbQ_{Y_U}) 
$
since for any smooth proper family of complex varieties the direct images define variations of Hodge structures~\cite[chapt.~III]{VoisinHodge1}. To construct the subvariations $\bbV_\pm$ we consider the embedding into the constant family
\[
\xymatrix@C=3em@M=0.5em{
 Y_U \ar@{^{(}->}[r]^-i \ar[dr]_-{f_U} & U\times X \ar[d]^-{p_U} \\
 & U
}
\]
where $p_U$ is the projection onto the first factor and where the closed immersion~$i$ is defined by $i(x,y)=(x+y,y)$, recalling that $Y_U \subset \Theta \times \Theta$ by construction. For each $n$ we have by adjunction a morphism
\[
 H^n(X, \bbQ) \otimes_\bbQ \bbQ_U \;=\; R^np_{U*}(\bbQ_{U\times X}) \; \longrightarrow \; R^nf_{U*}(\bbQ_{Y_U})
\]
of polarized variations of pure Hodge structures which induces on stalks the usual restriction morphism. The involution
$
 \sigma:  U\times X  \rightarrow U\times X,  (x, y) \mapsto (x, x-y)
$
preserves $Y_U$ and restricts on each fibre $Y_x$ to the involution $\sigma_x$. So we can mimic the constructions of section~\ref{sec:intersections} to get $\bbV_\pm \subset \bbH^{g-2}$ with the desired property. \qed

\medskip

It is now natural to ask for the monodromy groups of the local systems underlying the variations $\bbV_\pm$ of Hodge structure. Recall that by definition these monodromy groups are the image of the monodromy representation of $\pi_1(U, x)$ on the stalks of~$\bbV_\pm$ at a chosen base point $x\in U(\bbC)$. The rest of this paper is concerned with the proof of theorem~\ref{thm:generic_monodromy} which can be rephrased in the following form:

\medskip

{\bf Theorem. }{\em
For a general complex ppav $X$ of dimension $g=4$, the monodromy group $G$ of $\bbV_+$ is a subgroup of index at most two in the Weyl group $W(E_6)$.}

\medskip

{\em Plan of the proof.} The lattice computations in section~\ref{sec:lattice} show that if we view~$\bbV_+$ as a polarized variation of $\bbZ$-Hodge structures, then its stalks can be identified with the lattice $E_6$ up to a rescaling. Hence
\[
 G \; \leq \; W(E_6)
\]
is a subgroup of the Weyl group, see proposition~\ref{prop:upper-bound}. If equality holds, then we are done. If not, then $G$ is contained in a maximal proper subgroup $M < W(E_6)$.~Now the table of maximal subgroups in~\cite[p.~26]{ATLAS} shows that any such subgroup $M$ must be conjugate to one of the following subgroups:
\smallskip
\begin{enumerate}
\item[\em (a)] the simple subgroup $W^+(E_6) = \ker(W(E_6) \stackrel{\sgn}{\longrightarrow} \{ \pm 1\} )$ of index two,
\vspace*{0.5em}
\item[\em (b)]
the stabilizer of a line through a root vector in the lattice $E_6$ or through a minimal vector in the dual lattice, 
\vspace*{0.5em}
\item[\em (c)] three other subgroups, two of order $2^4\cdot 3^4$ and one of order $2^7\cdot 3^9$. 
\end{enumerate}
\smallskip

\noindent
Any proper subgroup of the maximal subgroup in {\em (a)} is by loc.~cit.~also contained in some of the other maximal subgroups. So for the proof of theorem~\ref{thm:generic_monodromy} it will be enough to show that for a general ppav the monodromy group $G$ is not contained in any of the maximal subgroups in {\em (b)} and {\em (c)}.

\medskip

One tool for this will be the Tannakian formalism for the convolution of perverse sheaves developed in~\cite{KrWVanishing}. With notations as in loc.~cit.~let $\delta_\Theta = \bbC_\Theta[g-1]$ be the constant perverse sheaf supported on the theta divisor. We will see in lemma~\ref{lem:L_plus_delta} that over the open dense subset $U\subset X$ the local system underlying $\bbV_+$ can be identified with a direct summand of the convolution $\delta_\Theta * \delta_\Theta$. It is shown in~\cite{KrWSchottky} that for a general ppav of dimension $g=4$ the decomposition of such convolutions is controlled by the representation theory of the Tannaka group $G(\delta_\Theta) = \Sp_{24}(\bbC)$, and this will imply that the local system underlying $\bbV_+$ is irreducible. In other words the monodromy group~$G$ acts irreducibly on the stalks $E_6\otimes_\bbZ \bbC$ of this local system, see proposition~\ref{cor:irreducible-monodromy}. Hence the group $G$ cannot be contained in the stabilizer of a line as in case {\em (b)} of the above list. 

\medskip

The maximal subgroups in~{\em (c)} cannot be excluded with the same argument since some of them act irreducibly on $E_6\otimes_\bbZ \bbC$. However, a look at the group orders shows that none of these subgroups contains the alternating group $\frakA_6$, so theorem~\ref{thm:generic_monodromy} will follow if we can construct an embedding
\[
 \frakA_6 \; \hookrightarrow \; G.
\]
For this we use a degeneration argument to fill in the dotted lines in the following Hasse diagram, where $G(d_2)$ denotes the Galois group of the difference morphism for a general curve of genus four as in theorem~\ref{thm:general-difference} in the introduction.
\smallskip
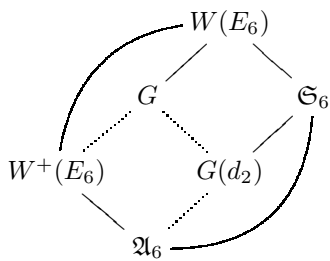
\[
\xymatrix@C=0.4em@R=1.3em@M=0.3em{
&& W(E_6) & \\
& G \ar@{-}[ur] && \frakS_6 \ar@{-}[ul]\\
W^+(E_6) \ar@{..}[ur] \ar@/^1.7pc/@{-}[uurr] && G(d_2) \ar@{..}[ul] \ar@{-}[ur] & \\
& \frakA_6 \ar@{-}[ul] \ar@{..}[ur]  \ar@/_1.85pc/@{-}[uurr] &&
}
\]
More specifically, as in~\cite[sect.~5]{KrWSchottky} one finds over a smooth complex quasi-projective curve $S$ a principally polarized abelian scheme $X_S \rightarrow S$ equipped with a relative theta divisor $\Theta_S \subset X_S$ such that the following two properties hold.

\smallskip
\begin{itemize}
 \item For some $s_0 \in S(\bbC)$ the fibre $X_{s_0}$ is the Jacobian variety of a general curve of genus four, whereas for all $s\neq s_0$ the theta divisor $\Theta_s\subset X_s$ is smooth.
\medskip
 \item The relative addition morphism $f: \; Y_S = \Theta_S \times_S \Theta_S \longrightarrow X_S$ restricts to a smooth proper morphism 
\[
 f_{U_S}: \quad f^{-1}(U_S) \; \longrightarrow \; U_S
\]
over some Zariski-open dense subset $U_S \subset X_S$ that surjects onto $S$.
\end{itemize}
\smallskip
The second property is the relative version of the Bertini-type lemma~\ref{lem:bertini} and can be achieved via the generic flatness theorem~\cite[th.~6.9.1]{EGA-IV-2} and the fibrewise flatness criterion~\cite[th.~11.3.10]{EGA-IV-3}, using the fact that for any proper flat morphism of varieties and any point of the target over which the fibre is smooth, there exists an open neighborhood over which the morphism is smooth~\cite[ex.~III.10.2]{Hartshorne}. In this relative context, over the open dense subset $U_S\subset X_S$ we have variations  of Hodge structures $\bbV_{S\pm}$ which on the fibre $U_s \subset X_s$ over each point $s\in S(\bbC)$ restrict to~$\bbV_\pm$ as defined in lemma~\ref{lem:vhs}. If we denote by $G(s)$ the monodromy group of the restriction $(\bbV_{S+})|_{U_s}$ to such a fibre, then a general result about degenerating monodromy representations (recalled in lemma~\ref{lem:degenerate_monodromy}) says that for general $s\in S(\bbC)$ we have an embedding
\[
 G(s_0) \; \hookrightarrow \; G(s).
\]
So we will be finished if we can show that the monodromy group $G(s_0)$ for the Jacobian variety $X_{s_0}=JC$ of a general curve $C$ of genus $g=4$ contains the group~$\frakA_6$ as a subgroup. To achieve this we show in proposition~\ref{prop:jacobian} that $G(s_0)$ is the Galois group of the difference morphism
\[
 d_2: \quad C_2 \times C_2 \; \longrightarrow \; \Pic^0(C), \quad (D, E) \mapsto [ \calO_C(D-E) ]
\]
where $C_2$ denotes the second symmetric product of the curve. In section~\ref{sec:diff} we will see that the Galois group $G(d_2)$ of this difference morphism is either the alternating group $\frakA_6$ or the full symmetric group $\frakS_6$ as claimed in theorem~\ref{thm:general-difference}. This fills in the dotted lines in the above Hasse diagram. \qed

\medskip

\section{The upper bound $W(E_6)$} \label{sec:upper-bound}

Let $X$ be a complex ppav of dimension $g=4$ with a smooth theta divisor, and with notations as in lemma~\ref{lem:vhs} consider the variation of $\bbQ$-Hodge structures $\bbV_+$ on the open dense  $U\subset X$. Repeating the above constructions  on the level of integral rather than rational cohomology, one sees that $\bbV_+$ underlies a polarized variation of~$\bbZ$-Hodge structures in a natural way. Proposition~\ref{prop:E6} identifies the fibres of this variation of Hodge structures with the lattice $K^\perp \cong E_6(-1)$ up to a rescaling. 

\medskip

Let $G$ be the monodromy group attached to $\bbV_+$ as above. Since the monodromy action for any smooth proper family of varieties preserves the intersection form on the cohomology of the fibres, $G$ is contained in the automorphism group $\Aut(E_6)$ of the $E_6$-lattice. Recall~\cite[p.~126]{CS2} that we have a product decomposition
\[
 \Aut(E_6) \;=\; W(E_6) \times \{\pm 1\}
\]
where the subgroup $\{ \pm 1\}$ acts via multiplication by $\pm 1$ on the $E_6$-lattice.

\begin{prop} \label{prop:upper-bound}
With notations as above, the monodromy group~$G$ is already contained in the subgroup \[ W(E_6) \; < \; \Aut(E_6). \]
\end{prop}

{\em Proof.} Fix a base point $x\in U(\bbC)$, and consider the surface $Y=Y_x$. With the notations of section~\ref{sec:lattice} we are interested in the monodromy action of~$\pi_1(U, x)$ on the lattice $H^2(Y, \bbZ)^+ \cong L(2)$. On the sublattice $H^2(X, \bbZ) \cong K(2)$ this action is trivial. So after dividing the bilinear forms by two, we obtain from proposition~\ref{prop:summary} an embedding
\[
 M\;=\; K \oplus E_6(-1)  \; \hookrightarrow \; L\;=\; U^{13} \oplus E_8(-1)
\]
with an action of $\pi_1(U, x)$ on $L$ that is trivial on $K$. On the orthocomplement $E_6(-1)=K^\perp$ this defines an epimorphism 
$\pi_1(U, x) \twoheadrightarrow  G \hookrightarrow  \Aut(E_6)$.
To show that the image $G$ of this epimorphism is already contained in the subgroup $W(E_6)$, we consider the action of $G$ on the discriminant group of our lattice.

\medskip

Since the lattices $L$ and $M$ both have the same rank, $L$ is obtained from its sublattice $M$ by adjoining certain glue vectors, by which we mean as in~\cite{CS2} vectors from the dual lattice 
\[
 M^* \;=\;
 \bigl \{ m\in M\otimes_\bbZ \bbQ \mid b( m, n) \in \bbZ \; \; \forall n \in M \bigr \}
 \;=\; K^* \oplus E_6(-1)^* 
\]
where $b$ denotes the bilinear form of the lattice $M$ extended to $M\otimes_\bbZ \bbQ$. Now the lattice~$E_6(-1)$ has discriminant $3$ whereas $L$ is unimodular, so there exists inside the sublattice $L \subseteq M^*=K^* \oplus E_6(-1)^*$ at least one glue vector of the form 
\[
 \lambda \;=\; k + e \quad \textnormal{with} \quad  k\in K^*, \; e\in  E_6(-1)^* \quad \textnormal{but} \quad e \; \notin \; E_6(-1).
\]
Fixing $k$ as above, put
$S = \{ f \in E_6(-1)^* \mid k+f \in L \}$.
This set $S$ is stable under the monodromy operation, indeed this operation fixes $k$ because $H^2(X, \bbZ) \cong K(2)$ is contained in the monodromy invariant part of the cohomology. By construction~$S$ contains the coset $e + E_6(-1)$. On the other hand, for any $f\in S$ we have 
\[ f-e \;=\; 
 (k+f)-(k+e)
 \;\in\; E_6(-1)^* \cap L 
 \;=\; E_6(-1) 
\] 
where the last equality uses that $E_6(-1) = K^\perp$ is primitive in $L$ because it is the orthocomplement of a sublattice. Altogether then
$S = e + E_6(-1)$
is a non-trivial coset of $E_6(-1)$ in $E_6(-1)^*$ which is preserved by the monodromy operation. Since the discriminant group
\[
 E_6(-1)^*/E_6(-1) \;\cong\; E_6^*/E_6 \;\cong\; \bbZ/3\bbZ
\]
is generated by any non-trivial coset, it follows that $G$ operates trivially on this discriminant group. Now recall~\cite[p.~27]{ATLAS} that inside~$\Aut(E_6)=W(E_6) \times \{\pm 1\}$ we have
$
 W(E_6) = \ker (\Aut(E_6) \rightarrow \Aut(E_6^*/E_6) ),
$
so the fact that the monodromy group~$G$ acts trivially on the discriminant group implies that $G\leq W(E_6)$. 
\qed

\bigskip

Let us briefly explain the connection with the $27$ {\em Prym-embedded curves} that have been studied by E.~Izadi~\cite{Iz}, though this will not be used in the sequel. If $X$ is a general complex ppav of dimension $g=4$, then by~\cite[prop.~6.4]{BeTnull} it is the Prym variety 
$
 X \cong \mathit{Prym}(\tilde{C}/C)
$
of an \'etale double cover of smooth projective curves of genus $9$ resp.~$5$. Consider the covering involution $\iota: \tilde{C} \rightarrow \tilde{C}$, and let $\alpha: \tilde{C} \hookrightarrow J\tilde{C}$ be a translate of the Abel-Jacobi map. The morphism
\[
 \tilde{C} \; \longrightarrow \; J\tilde{C}, \quad p \; \mapsto \; \alpha(p)-\alpha(\iota(p))
\]
factors over the kernel of the norm map $J\tilde{C} \rightarrow JC$, so the image of some translate of this morphism lies in the connected component $\mathit{Prym}(\tilde{C}/C)$ of this kernel. Via the chosen isomorphism $X\cong \mathit{Prym}(\tilde{C}/C)$ this defines an embedding $\tilde{C} \hookrightarrow X$ which up to a translation is determined uniquely by the \'etale double cover. Embeddings of this form are known as Abel-Prym embeddings, and by a Prym-embedded curve in~$X$ we mean the image of any such embedding. Note that the covering $\tilde{C} \rightarrow C$ is not determined  uniquely by $X$, indeed a dimension count shows that the general fibre of the Prym morphism
$
 \calR_5 \rightarrow \calA_4
$
has dimension two. In this context E. Izadi has observed~\cite[cor.~4.9]{Iz} that for general $x\in X(\bbC)$ precisely $27$ Prym-embedded curves are contained in $Y=Y_x = \Theta \cap \Theta_x$. This number suggests a relationship with the lattice $E_6$ in the cohomology of this surface. Recall from~\cite{CS2} that the dual lattice $E_6^*$ has $27$ pairs of minimal vectors with norm $4/3$. It turns out that the fundamental classes of Prym-embedded curves define such minimal vectors and are related to the glue vector $\lambda = k + e \in L$ with components $k\in K^*$ and $e\in E_6(-1)^*$ from the proof of proposition~\ref{prop:upper-bound} in the following way.

\begin{lem}
The above glue vector $\lambda = k + e$ can be taken to be the fundamental class 
$[ \tilde{C}  ] \in  H^2(Y, \bbZ)^+  \cong  L(2) $
of a Prym-embedded curve $\tilde{C}$, and then $e\in E_6(-1)^*$ is a minimal vector which has norm $-4/3$.
\end{lem}

{\em Proof.} Any integral cohomology class $\lambda \in H^2(Y, \bbZ)^+$ can be written uniquely as a sum of rational classes
\[ 
 \lambda \;=\; \alpha + \beta \quad \textnormal{with} \quad \alpha \;\in\; H^2(X, \bbQ) \quad \textnormal{and} \quad \beta \;\in\; H^2(X, \bbQ)^\perp \subset H^2(Y, \bbQ)^+
\]
and it defines a glue vector iff $\alpha \notin H^2(X, \bbZ)$. The integrality of~$\alpha$ can be checked via the Gysin morphism for the embedding $i: Y\hookrightarrow X$ since $i_!(\lambda) = \alpha \cup [\Theta ]^2 \in H^6(X, \bbQ)$ is the image of $\alpha$ under the Lefschetz isomorphism in the following diagram.
\smallskip
\[
\xymatrix@M=0.5em@R=0.8em@C=1.5em{
 H^2(X, \bbQ) \ar[rr]^-\cong \ar[dr]^-{i^*} && H^6(X, \bbQ) \\
 & H^2(Y, \bbQ)^+ \ar[ur]^-{i_!} & \\
 H^2(X, \bbZ) \ar@{-}[r] \ar[dr] \ar@{^{(}->}[uu] & \ar[r] & H^6(X, \bbZ) \ar@{^{(}->}[uu] \\
 & H^2(Y, \bbZ)^+ \ar[ur] \ar@{^{(}->}[uu] &
}
\]
By~\cite[sect.~12.2]{BL} the fundamental class $\lambda = \alpha + \beta \in H^2(Y, \bbQ)^+$ of a Prym-embedded curve satisfies
\[ \alpha \cup [\Theta]^2 \;=\; i_!(\lambda) \;=\; [\Theta]^3 /3 \;\notin\; H^6(X, \bbZ) \]
and hence indeed provides a glue vector for our lattices. It also follows from the above that
$
 \deg_Y(\alpha^2) = \deg_X(\alpha^2\cup [\Theta]^2) = \deg_X([\Theta]^4/9) = 8/3
$
where the last equality is due to the Poincar\'e formula. Furthermore, by definition $\lambda \in H^2(Y, \bbQ)$ is the class of a curve $\tilde{C}\subset Y$ of genus $g_{\tilde{C}}=2g+1=9$. Hence the adjunction formula for the genus of a smooth curve on a smooth projective surface~\cite[prop.~V.1.5]{Hartshorne} gives the intersection number
\[
 \deg_Y(\lambda^2) \;=\; 2g_{\tilde{C}} - 2 - \deg_Y(K_Y|_{\tilde{C}}) \;=\; 16 - 16 \;=\; 0,
\]
indeed for the canonical class $K_Y$ lemma~\ref{lem:chern-classes} and the formula $i_!(\lambda)=[\Theta]^3/3$ from above imply
$\deg_{\tilde{C}}(K_Y|_{\tilde{C}}) = \deg_Y(\lambda \cup K_Y)= \deg_X(2\cdot [\Theta]^4 / 3) = 16$.
Finally we have $\deg_Y(\beta^2) = \deg_Y(\lambda^2)-\deg_Y(\alpha^2) = -8/3$, and to see that the vector $e$ is minimal it only remains to divide the bilinear form by two. \qed

\medskip

It seems tempting to prove theorem~\ref{thm:generic_monodromy} directly by looking at the monodromy operation on the fundamental classes of the Prym-embedded curves, inspired by the result of~\cite{Do}. In what follows we give a different proof, and from this point of view the theorem can be seen as an independent result about the monodromy of the Prym-embedded curves.

\medskip

\section{Negligible constituents} \label{sec:irreducible}

We now show that for a general ppav the local systems underlying the variations of Hodge structures $\bbV_\pm$ from lemma~\ref{lem:vhs} are irreducible. In fact our results are stronger since we consider perverse sheaves on the whole ppav rather than local systems on a open dense subset. This shifted perspective reveals a close connection with the Tannakian formalism of~\cite{KrWVanishing} which may be of independent interest.

\medskip

For any complex abelian variety $X$ the group law $a: X\times X\to X$ defines a {\em convolution product} on the derived category $\Dbc(X, \bbC)$ of bounded constructible sheaf complexes in the sense of~\cite{BBD} by
\[
 K*L \;=\; Ra_*(K\boxtimes L) \;\in \; \Dbc(X, \bbC) \quad \textnormal{for} \quad K, L\in \Dbc(X, \bbC).
\]
Let $\Perv(X, \bbC) \subset \Dbc(X, \bbC)$ be the full subcategory of perverse sheaves. In~\cite{KrWVanishing} it has been shown that for all $n\neq 0$ and all $P, Q\in \Perv(X, \bbC)$ the perverse cohomology sheaves 
$\pH^n(X, P*Q)\in \Perv(X, \bbC)$
are {\em negligible} in the sense that all their subquotients have Euler characteristic zero. By way of contrast, let us say that a perverse sheaf is {\em clean} if it does not have any negligible subquotients. 

\medskip

Suppose that $X$ has dimension $g\geq 2$ and is principally polarized. Let~$i: \Theta\hookrightarrow X$ be the inclusion of a symmetric theta divisor, $j: \Theta^{sm} \hookrightarrow \Theta$ the open immersion of its smooth locus, and consider the perverse intersection cohomology sheaf 
\[ \delta_\Theta \;=\; \IC_\Theta[g-1] \;=\; i_* j_{!*} (\bbC_{\Theta^{sm}}[g-1]) \;\in\; \Perv(X, \bbC) \]
which is the intermediate extension of the constant sheaf on the smooth locus.

\begin{lem} \label{lem:convolution-with-theta}
Suppose that the theta divisor $\Theta$ has at most isolated singularities and that $P\in \Perv(X, \bbC)$ is clean. Then for all $n\in \bbZ$ any simple negligible subquotient of the perverse cohomology sheaf $\pH^n(\delta_\Theta * P)$ is isomorphic to $\delta_X = \bbC_X[g]$.
\end{lem}

{\em Proof.} By Verdier duality we can assume $n\leq 0$. Since the theta divisor is a local complete intersection,  by~\cite[III.6.5]{KW} the shifted constant sheaf $\lambda_\Theta = \bbC_\Theta[g-1]$ is perverse, so we have an exact sequence of perverse sheaves
\[
 0 \longrightarrow \kappa \longrightarrow \lambda_\Theta \longrightarrow \delta_\Theta \longrightarrow 0
\]
where $\kappa$ is a skyscraper sheaf supported on the finite set of singular points of the theta divisor. Convolution with the clean perverse sheaf $P$ yields a distinguished triangle
$\kappa * P \to
 \lambda_\Theta * P \to
 \delta_\Theta * P \to  \cdots
$
where $\kappa * P$ is clean and perverse. By the long exact sequence of perverse cohomology groups, it will therefore suffice to show that every simple negligible subquotient of $\pH^n(\lambda_\Theta*P)$ for $n\leq 0$ is isomorphic to $\delta_X$. But this follows as in \cite[ex.~3]{WeT} from the excision sequence for the ample divisor $\Theta \subset X$ with affine complement $X\setminus \Theta$, using Artin's vanishing theorem. \qed

\medskip

Now consider the convolution $\delta_\Theta * \delta_\Theta$. By the decomposition theorem this is a direct sum of degree shifts of semisimple perverse sheaves. So we can write
\[
 \delta_\Theta * \delta_\Theta \;=\; S^2(\delta_\Theta) \oplus \Lambda^2(\delta_\Theta)
\]
where the symmetric and the alternating square on the right are the maximal direct summands on which the commutativity constraint 
$ 
  S = S_{\delta_\Theta, \delta_\Theta}: \delta_\Theta * \delta_\Theta  \rightarrow  \delta_\Theta * \delta_\Theta
$
from \cite{WeBN} acts by $+1$ resp.~by $-1$. We want to relate these sheaf complexes to the variations of Hodge structures $\bbV_\pm$ in lemma~\ref{lem:vhs}. In this context the contribution from the weak Lefschetz theorem will be described by the complexes
\[
 \bbL_\pm \; \;=\;  \; \bigoplus_{\; \mu \, \in \, I_\pm} H^{g-2-|\mu|}(X, \bbC) \otimes_\bbC \delta_X[\mu]
 \; \; \; \; \textnormal{with} \; \; \;
 \begin{cases}
  I_+ = 1+2\bbZ, \\
  I_- = 2\bbZ,
 \end{cases}
\]
where by convention $H^i(X, \bbC)=0$ for $i<0$. Strictly speaking certain Tate twists should be inserted in the above definition, but we will ignore these.

\begin{lem} \label{lem:L_plus_delta}
If the theta divisor $\Theta$ has at most isolated singularities, then there are decompositions
\[
 S^2(\delta_\Theta) \;\cong\; \bbL_+ \oplus \delta_+ 
 \quad \textnormal{\em and} \quad
 \Lambda^2(\delta_\Theta) \;\cong\; \bbL_- \oplus \delta_-
\]
where $\delta_\pm \in \Perv(X, \bbC)$ are semisimple perverse sheaves, and putting $\epsilon = (-1)^{g-1}$ we have isomorphisms
$
\delta_{\pm}|_U  \cong \bbV_{\pm \epsilon}[g ]
$ 
over the open dense $U\subset X$ in lemma~\ref{lem:vhs}.
\end{lem}

{\em Proof.} The group law of the ppav induces a morphism $f: \Theta \times \Theta \to X$ whose restriction over $U\subset X$ is the family $f_U: Y_U \rightarrow U$ in lemma~\ref{lem:bertini}. Since the theta divisor has at most isolated singularities, we can assume that
$f^{-1}(U) \subseteq \Theta^{sm} \times \Theta^{sm}$
where $\Theta^{sm} \subseteq \Theta$ denotes the smooth locus of the theta divisor. The definitions then imply that
$(\delta_\Theta * \delta_\Theta)|_U \cong Rf_{U*}\bbC_{Y_U}[2g-2]$.
By the decomposition theorem this direct image is a direct sum of degree shifts of semisimple perverse sheaves and hence equal to $\bigoplus_{n \in \bbZ} \, \bbH^n[2g-2-n]$ for the local systems $\bbH^n = R^nf_{U*}(\bbC_{Y_U})$. As in lemma~\ref{lem:vhs} we then get
\[ 
 (\delta_\Theta * \delta_\Theta)|_U \;\cong\; (\bbL_+ \oplus \bbL_-)|_U \oplus \bbV[g]
 \quad \textnormal{where} \quad \bbV \;=\; \bbV_+ \oplus \bbV_-. \]
The decomposition theorem thus implies $\delta_\Theta * \delta_\Theta \cong \bbL_+ \oplus \bbL_- \oplus \delta$ for some sheaf complex $\delta \in \Dbc(X, \bbC)$. But for~$n\neq 0$ the perverse cohomology sheaves $\pH^n(\delta_\Theta*\delta_\Theta)$ must be negligible and hence by lemma~\ref{lem:convolution-with-theta} multiples of $\delta_X$. Since over the open dense $U\subset X$ we have $\delta|_U \cong \bbV[g]$, it follows that $\delta$ is a perverse sheaf. 

\medskip

Now consider the decomposition $\delta = \delta_+\oplus \delta_-$ where $\delta_\pm \subset \delta$ are the maximal perverse subsheaves on which the commutativity constraint~$S$ acts by $\pm 1$. Going back to the definitions, one checks that this commutativity constraint induces on the stalks the involution
$\epsilon \cdot \sigma_x^*: H^{n}(Y_x, \bbQ) \rightarrow H^{n}(Y_x, \bbQ)$
where $\sigma_x: Y_x\rightarrow Y_x$ is the involution that we defined in section~\ref{sec:translates} and where the twist by the sign $\epsilon = (-1)^{g-1}$ comes from the degree shift in $\delta_\Theta = \IC_\Theta[g-1]$ by the Koszul rule. \qed

\medskip

\begin{lem} \label{lem:negligible}
If the theta divisor $\Theta \subset X$ is smooth, then $\delta_\pm$ are clean. 
\end{lem}

{\em Proof.} If this were not true, then $\delta_+\oplus \delta_-$ would contain a direct summand $\delta_X$ by lemma~\ref{lem:convolution-with-theta}. Then the inclusion
$
 H^{-g}(X, \bbL_+ \oplus \bbL_-) \; \hookrightarrow \; H^{-g}(X, \delta_\Theta * \delta_\Theta)
$
would be strict. However, the K\"unneth formula and the weak Lefschetz theorem for the smooth ample divisor $\Theta \subset X$ show that
\begin{align} \nonumber
 H^{-g}(X, \delta_\Theta * \delta_\Theta) & \;\;=\;\; 
 \bigoplus_{n=1}^{g-1} \; H^{-n}(X, \delta_\Theta) \otimes H^{n-g}(X, \delta_\Theta) \\ \nonumber
 & \;\;=\;\; \bigoplus_{n=1}^{g-1} \; H^{g-1-n}(X, \bbC) \otimes H^{n-1}(X, \bbC)
\end{align}
which by direct inspection is equal to $H^{-g}(X, \bbL_+ \oplus \bbL_-)$. \qed

\medskip

For a general complex ppav of dimension $g$ the theta divisor is smooth, so the above lemma leads to the final result of this section.

\begin{prop} \label{cor:irreducible-monodromy}
For a general complex ppav $X$ of dimension $g$ the local systems underlying $\bbV_\pm$ are irreducible.
\end{prop}

{\em Proof.}  For a general ppav we know from~\cite{KrWSchottky} that the Tannaka group $G(\delta_\Theta)$ is either a symplectic or a special orthogonal group and that the perverse intersection cohomology sheaf $\delta_\Theta$ corresponds to the standard representation of this group. Now the symmetric and alternating square of the standard representation are irreducible up to a trivial one-dimensional representation, and the latter corresponds to the skyscraper sheaf $\delta_0$ of rank one supported in the origin $0\in X(\bbC)$. Furthermore, each representation of the Tannaka group corresponds to a unique clean perverse sheaf. Hence from lemma~\ref{lem:negligible} we obtain that there exist clean and simple perverse sheaves $\gamma_\pm \in \Perv(X, \bbC)$ such that $\delta_{+\epsilon} = \gamma_{+\epsilon} \oplus \delta_0$ and $\delta_{-\epsilon} = \gamma_{-\epsilon}$ for $\epsilon = (-1)^{g-1}$ and by lemma~\ref{lem:L_plus_delta} the simplicity of $\gamma_\pm$ in particular implies that the local systems underlying~$\bbV_\pm$ are irreducible. \qed

\medskip

\section{Jacobian varieties} \label{sec:bn}

Let us now see what happens if the complex ppav $X=JC$ is the Jacobian variety of a smooth projective curve $C$ of genus $g\geq 2$. For $n\in \bbN$ let us denote by~$C_n = (C\times \cdots \times C)/\frakS_n$ the $n$-fold symmetric product of the curve and by $W_n\subset X$ the image of~$C_n$ under a suitable translate of the Abel-Jacobi map, normalized such that the theta divisor $\Theta = W_{g-1}$ becomes symmetric. We write $f_n: C_n \rightarrow W_n$ for the Abel-Jacobi map and $g_n: W_n \times W_n \rightarrow X$ for the difference morphism given by~$g_n(x,y)=y-x$. Then the composite morphism
\[ 
d_n: \;\;
 \xymatrix@C=3.3em@M=0.6em{
 C_n \times C_n \ar[r]^-{(f_n, f_n)} 
 & W_n\times W_n \ar[r]^-{g_n}
 & X \; \cong \; \Pic^0(C)
}
\]
does not depend on the chosen normalization of the Abel-Jacobi map, and it sends a pair $(D, E)\in  C_n \times C_n$ of effective divisors of degree $n$ on the curve to the isomorphism class of the line bundle $\calO_C(D-E)$. 

\medskip

In what follows we will be especially interested in non-hyperelliptic curves of genus $g=4$. For these the theta divisor has only isolated singularities, so we can consider the variations of Hodge structures $\bbV_\pm$ from lemma~\ref{lem:vhs} and the perverse sheaves $\delta_\pm$ from lemma~\ref{lem:L_plus_delta}. Here we have the following situation.

\begin{prop} \label{prop:jacobian} \label{lem:jacobian-const}
For the Jacobian variety $X=JC$ of a non-hyperelliptic smooth projective curve $C$ of genus $g=4$ the following properties hold. 
\begin{enumerate}
\item[\em (a)] The semisimple perverse sheaves $\delta_\pm$ each contain $\delta_X$ precisely once. \smallskip
\item[\em (b)] After shrinking the open dense $U\subset X$ we have
$ \bbV_+  \cong  d_{2*}(\bbC_{C_2\times C_2})|_U$.
\end{enumerate}
\end{prop}

{\em Proof.} For {\em (a)} we must show $\dim (H^{-g}(X, \delta_+)) = \dim (H^{-g}(X, \delta_-)) = 1$. For this we can use the same computation as in lemma~\ref{lem:negligible}, the only difference is that for a Jacobian variety the hypercohomology $H^\bullet(X, \delta_\Theta)$ is larger than the one for a general ppav. Indeed it has been shown in~\cite[cor.~13(iii)]{WeBN} that if the curve $C$ is not hyperelliptic, then
\[
 H^{i}(X, \delta_\Theta) \;=\;  
 \begin{cases}
  \; H^{i+3}(X, \bbC) & \textnormal{for $i\in \{ -3, -2\}$}, \\
  \; H^{i+3}(X, \bbC) \oplus H^{i+1}(X, \bbC) & \textnormal{for $i \in \{-1, 0\}$. }
 \end{cases}
\]
Compared with lemma~\ref{lem:negligible}, for $i=-1$ the extra summand $H^0(X, \bbC) \subset H^{-1}(X, \delta_\Theta)$ gives a two-dimensional extra term
\[
 H^0(X, \bbC) \otimes H^{-3}(X, \delta_\Theta) 
 \; \oplus \; H^{-3}(X, \delta_\Theta) \otimes H^0(X, \bbC) 
 \; \subset \; H^{-4}(X, \delta_\Theta * \delta_\Theta).
\]
This extra term splits into two one-dimensional contributions, one in the symmetric and one in the alternating tensor square. Hence 
\begin{align} \nonumber
 H^{-4}(X, S^2(\delta_\Theta)) & \;=\;  H^{-4}(X, \bbL_+) \oplus \bbC, \\ \nonumber
 H^{-4}(X, \Lambda^2(\delta_\Theta)) & \;=\; H^{-4}(X, \bbL_-) \oplus \bbC,
\end{align}
and part {\em (a)} follows because $S^2(\delta_\Theta)=\bbL_+\oplus \delta_+$ and $\Lambda^2(\delta_\Theta)=\bbL_-\oplus \delta_-$. Part~{\em (b)} is most conveniently checked in the setting of Brill-Noether sheaves and will be postponed after the proof of lemma~\ref{lem:one-plus-irred} below. \qed

\medskip

Before we proceed further, let us recall some basic facts from~\cite{WeBN}. Let $X=JC$ be the Jacobian variety of a smooth non-hyperelliptic complex projective curve of genus $g\geq 2$, and $C\hookrightarrow X$ a translate of the Abel-Jacobi map such that the corresponding theta divisor $\Theta = W_{g-1}$ is symmetric. By loc.~cit.~the convolution powers of the perverse intersection cohomology sheaf $\delta_C = \bbC_C[1]$ are described by the representation theory of the Tannaka group
\[ 
 G(\delta_C) \;=\; \Sl_{2g-2}(\bbC). 
\]
For partitions $\alpha = (\alpha_1, \alpha_2, \dots)$ such that $2g-2\geq \alpha_1 \geq \alpha_2 \geq \cdots \geq 0$, one can define the {\em Brill-Noether sheaf} 
\[ \pdelta_\alpha \;\in\; \Perv(X, \bbC) \]
as the simple perverse sheaf which corresponds via the Tannakian formalism of~\cite{KrWVanishing} to the irreducible representation of~$\Sl_{2g-2}(\bbC)$ whose highest weight is given by the conjugate partition $\alpha^t$ in the basis of the weights defined by the diagonal entries of matrices. In particular, the singleton partitions $\alpha = (\alpha_1)$ correspond to the fundamental representations, and for these it has been shown in loc.~cit.~that the above construction leads to the perverse intersection cohomology sheaves
\[ 
  \pdelta_i \;=\; \delta_{\, W_i}
  \quad \textnormal{and} \quad
  \pdelta_{2g-2-i} \; = \; \delta_{-W_i}
  \quad \textnormal{for} \quad 0\leq i \leq g-1,
\]
where the second equality is due to the Riemann-Roch theorem. The decomposition of arbitrary convolution products between these perverse intersection cohomology sheaves can be obtained via the {\em Littlewood-Richardson rule} for the decomposition of tensor products of irreducible representations of $\Sl_{2g-2}(\bbC)$, see~\cite[sect.~9.3.5]{GW} and~\cite[sect.~5.4]{WeBN}. For the fundamental representations this gives 
\[
 \pdelta_m * \pdelta_n \;\; \cong \;\; \tau_{m,n} \; \oplus \; \bigoplus_{i=0}^n \; \pdelta_{(m+i, n-i)}
 \quad \textnormal{for} \quad
 m\geq n,
\]
where $\tau_{m,n} \in \Dbc(X, \bbC)$ is a negligible sheaf complex in the sense of the previous section. It has been shown in lemma~27 of loc.~cit.~that $\tau_{m,n}$ and more generally any negligible direct summand of a convolution of Brill-Noether sheaves is a direct sum of degree shifts of $\delta_X$. As an application we have the following result.

\begin{lem} \label{lem:one-plus-irred}
Let $X=JC$ be the Jacobian of a non-hyperelliptic smooth  projective curve of even genus $g=2n$. Then over some open dense $U\subset X$ there exists an irreducible local system $L_U$ of complex vector spaces such that
\[ d_{n*}(\bbC_{C_n\times C_n})|_U \cong \bbC_U \oplus L_U. \] 
\end{lem}

{\em Proof.} An argument like in lemma~\ref{lem:bertini} shows that over some open dense $U\subset X$ the morphism $d_n$ restricts to a finite \'etale cover. Furthermore, if $S$ denotes the singular locus of $W_n$, then the closed subset $(S\times W_n)\cup (W_n\times S)\subset W_n \times W_n$ is mapped under the difference morphism 
\[ g_n: \quad W_n\times W_n \; \longrightarrow \; X \] 
to a proper closed subset of $X$ by dimension reasons, so we can assume that $g_n^{-1}(U)$ is contained in $(W_n\setminus S)\times (W_n \setminus S)$. Now recall~\cite[p.~348]{GH} that by the singularity theorem of Riemann-Kempf the Abel-Jacobi morphism $f_n: C_n \rightarrow W_n$ restricts over the smooth locus to an isomorphism 
\[ f_n^{-1}(W_n\setminus S) \; \stackrel{\cong}{\longrightarrow} \; W_n\setminus S. \]
Putting everything together and using the commutative diagram
\[
\xymatrix@C=4em@R=2em@M=0.5em{
C_n \times C_n \ar[r]^-{(f_n,f_n)} \ar@/_1pc/[dr]^-{d_n} & W_n\times W_n \ar[d]_-{g_n} \ar@{^{(}->}[r]^-{(-i_n, i_n)} & X\times X  \ar@/^1pc/[dl]_-{a^{}} \\
& X &
}
\]
for the closed embedding $i_n: W_n\hookrightarrow X$, we get 
$
  d_{n*}(\bbC_{C_n\times C_n}[g])|_U \cong (\delta_{\, W_n}*\delta_{-W_n})|_U
$
by the definition of the convolution product. To control the right hand side, recall that $\delta_{\, W_n} = \pdelta_n$ and $\delta_{-W_n}=\pdelta_{3n-2}$ where the second equality is due to the Riemann-Roch theorem. Hence the Littlewood-Richardson rule says
\[
 \delta_{\, W_n} * \delta_{-W_n} \;\cong\; \pdelta_n * \pdelta_{3n-2} \;\cong\;
 \tau \; \oplus \; \bigoplus_{i=0}^n \; \pdelta_{3n-2+i, n-i} 
\]
where $\tau$ is a sum of degree shifts of $\delta_X$. In fact $\tau = \delta_X$ by adjunction since~$d_n$ is generically finite. So the claim of the lemma will follow for the cohomology sheaf~$ L_U = \calH^{-g}(\pdelta_{3n-2, n})|_U$ provided that for all $i>0$ the $\pdelta_{3n-2+i, n-i}$ have their support inside a proper closed subset of $X$. But this is indeed the case because by the Littlewood-Richardson rule the perverse sheaf
\[
 \bigoplus_{i=1}^n \; \pdelta_{3n-2+i, n-i} \;\cong \; \pdelta_{n-1}*\pdelta_{2g-(n-1)} \;\cong\; \delta_{W_{n-1}}*\delta_{-W_{n-1}}
\]
has the support $d_{n-1}(C_{n-1}\times C_{n-1}) \subset X$ of dimension $2(n-1)<g$. \qed

\medskip

{\em Proof of part (b) in proposition~\ref{prop:jacobian}}. 
Again we use Brill-Noether sheaves. The perverse sheaf $\delta_\Theta = \delta_{W_3} = \pdelta_3$ corresponds by definition to the third fundamental representation of the Tannaka group $\Sl_6(\bbC)$. Representation theory shows that the alternating square of this representation decomposes into two irreducible pieces of highest weight $\alpha^t$ for $\alpha = (4,2)$ and $\alpha = (6)$. Hence we get
$\delta_- \cong \pdelta_{4,2} \oplus \pdelta_6 \oplus \tau$
where $\tau \in \Perv(X, \bbC)$ is some negligible perverse sheaf. In fact $\tau=\delta_X$ by part {\em (a)} of the proposition. Furthermore $\pdelta_6 = \delta_0$ is a skyscraper sheaf, hence for $U\subset X$ sufficiently small we have
$
 \bbV_+ \;\cong\; \calH^{-4}(\delta_-)|_U \;\cong\; \calH^{-4}(\pdelta_{4,2} \oplus \delta_X)|_U \;\cong\; d_{2*}(\bbC_{C_2\times C_2})
$
where the first isomorphism is due to lemma~\ref{lem:L_plus_delta} and where the last one comes from the proof of lemma~\ref{lem:one-plus-irred} for $n=2$.
\qed

\medskip

\section{The difference morphism} \label{sec:diff}

This section is logically independent from the rest of the paper and does not use Brill-Noether sheaves. As above let $C$ be a smooth complex projective curve of genus $g=2n$ for some $n\in \bbN$, and put
$C_n = (C\times \cdots \times C)/\frakS_n$.
Motivated by the previous section, our goal is to understand the properties of the difference morphism $d_n: C_n\times C_n\rightarrow X = JC$ from above. To begin with, we claim that this morphism is generically finite of degree
\[
 N \;=\; \deg(d_n) \;=\; \tbinom{2n}{n} \, . 
\]
Indeed, by birationality of the Abel-Jacobi map $C_n\to W_n$ it suffices to check that the difference morphism $g_n: W_n\times W_n \to X$ is generically finite of degree $N$. The fibre of $g_n$ over a point $x\in X(\bbC)$ is isomorphic to $Z=W_n\cap (W_n+x)$. The Poincar\'e formula \cite[sect.~11.2.1]{BL} for the classes 
\[ [W_n] \;=\; [W_n + x] \;=\; \tfrac{1}{n!} \cdot [\Theta]^n \;\in\; H^n(X, \bbZ) \]
shows that for a sufficiently general point $x\in X(\bbC)$ the intersection $Z$ is finite of cardinality $N$, and our claim follows.

\medskip

Now pick $x\in U(\bbC)$, and consider the monodromy action of $\pi_1(U, x)$ on the $N$ distinct points of the fibre $d_n^{-1}(x)$. Labelling these points in any chosen order, we get a homomorphism
\[
 \pi_1(U, x) \; \longrightarrow \; \frakS_N
\]
to the symmetric group, and we define the Galois group $G(d_n)$ to be the image of this homomorphism. A different choice of the labelling only changes this subgroup by an inner automorphism of $\frakS_n$.

\begin{lem} \label{lem:hyp}
If $C$ is a hyperelliptic curve of genus $g=2n$, then the group $G(d_n)$ is isomorphic to the symmetric group~$\frakS_{2n}$, and its action on the $N$ points of a general fibre of $d_n$ can be identified with the natural permutation action of~$\frakS_{2n}$ on the set of $n$-element subsets of $\{1, 2, \dots, 2n\}$.
\end{lem}

{\em Proof.} If $g_2^1$ denotes the hyperelliptic linear series on $C$, then for every effective divisor $D\in C_n$ the linear series $n\cdot g_2^1 - D$ contains an effective divisor. Hence~$W_n$ is a translate of its negative $-W_n$, and it follows that up to a translation the difference morphism $d_n$ coincides with the addition morphism $C_n\times C_n \to X$. The latter factors as $C_n\times C_n  \to C_{2n} \to X$ where the Abel-Jacobi morphism $C_{2n} \to X$ is birational. So we must determine the monodromy group $H$ of the finite branched cover $C_n\times C_n \to C_{2n}$. This cover is not Galois, but obviously it is dominated by the Galois cover with group $\frakS_{2n}$ in the following commutative diagram.
\[
\xymatrix@M=0.6em@R=1.5em{
 C^{2n} \ar@{=}[r] \ar@/_1.5em/[ddr]_-{\frakS_{2n}} & C^n\times C^n \ar[d]^-{\; \frakS_n\times \frakS_n}  \\
  & C_n\times C_n \ar[d] \\
 & C_{2n}
}
\]
Take $p_1+\cdots + p_{2n}\in C_{2n}$ with all $p_i\in C$ distinct. The projection $C_n\times C_n \to C_n$ onto the first factor identifies the fibre $F$ of the cover $C_n\times C_n \to C_{2n}$ over this point with the set of all $n$-element subsets of~$\{ p_1, \dots, p_{2n}\}$, and the monodromy group $H$ is the image of the homomorphism
\[
 \varphi: \; \frakS_{2n} \longrightarrow \Aut(F)
\]
which is given by the action of the Galois group $\frakS_{2n}$ of $C^{2n}\to C_{2n}$ on the fibre~$F$ via permutation of  $p_1, \dots, p_{2n}$. Now $\varphi$ is injective, indeed the identity is the only permutation in $\frakS_{2n}$ that fixes all $n$-element subsets of~$\{ p_1, \dots, p_{2n}\}$. 
\qed

\medskip

Based on the hyperelliptic case we can now also deal with a general curve of even genus, proving theorem~\ref{thm:general-difference} from the introduction:

\medskip

{\bf Theorem. }{\em If $C$ is a general curve of genus $g=2n$, then the Galois group $G(d_n)$ is either the alternating group $\frakA_N$ or the full symmetric group $\frakS_N$.
}

\medskip

{\em Proof.}  The group $G(d_n)$ can be embedded into the symmetric group $\frakS_N$ in such a way that the monodromy representation is the restriction of the natural permutation representation $V=\bbC^N$ of the symmetric group. Note that $V$ splits as the direct sum of the one-dimensional trivial representation plus an irreducible representation $W$ of the symmetric group. By lemma~\ref{lem:hyp} we can assume $g>2$ so that a general curve of genus $g$ is non-hyperelliptic, and then lemma~\ref{lem:one-plus-irred} implies that the restriction of~$W$ to the subgroup $G(d_n) \leq \frakS_N$ remains irreducible.

\medskip

From the theory of permutation groups we know that any subgroup of $\frakS_N$ with this irreducibility property is $2$-transitive~\cite[th.~1(b)]{Saxl} and hence by~\cite[prop.~5.2]{Ca} contains a unique minimal normal subgroup $H$ which is either elementary abelian or simple. The elementary abelian case is ruled out by theorem~4.1 of loc.~cit.~since~$N$ is not a prime power. So by loc.~cit.~we are left with finitely many possible simple subgroups $H$ and we must see that in our case $H=\frakA_N$. By a direct inspection of the finitely many cases~\cite[appendix D]{KrDiss} this task can be reduced to the claim that for a general curve of genus $g=2n$ the monodromy group $G(d_n)$ contains the one from lemma~\ref{lem:hyp}. For this we use a degeneration into a hyperelliptic curve.

\medskip

Let $p: \calC \rightarrow S$ be a flat, projective family of smooth curves of genus~$g$ over a smooth quasi-projective complex curve $S$, and assume that the fibre $\calC_s = p^{-1}(s)$ is hyperelliptic for some point $s=s_0 \in S(\bbC)$ but non-hyperelliptic for all~$s\neq s_0$. The existence of such families follows e.g.~from~\cite[th.~XII.9.1]{ACG}. Let~$\calX$ be the relative Picard scheme of $\calC$ over $S$ as in~\cite[sect.~8.2]{BLR}, and $\calC_n = (\calC\times_S \cdots \times_S \calC)/\frakS_n$ the $n$-th relative symmetric product. We have a morphism $d_n: \calC_n \times_S \calC_n \to \calX$ which on the fibres over each $s\in S(\bbC)$ restricts to the difference morphism that we considered above. Let $U\subseteq \calX$ be an open subset over which $d_n$ is finite \'etale of degree $N$ and such that over every $s\in S(\bbC)$ the fibre $U_s = U\cap \calX_s$ is dense in $\calX_s$. The following lemma~\ref{lem:degenerate_monodromy} applies to $\bbL = d_{n*}(\bbC_{\calC_n\times_S \calC_n})|_U$ and shows that the Galois group of $d_n$ for a general curve contains the one for a hyperelliptic curve. \qed

\medskip

For the above degeneration argument we have used the following semicontinuity property for monodromy groups in families. Let $f: U\rightarrow S$ be a smooth morphism of complex algebraic varieties whose target is a quasi-projective curve $S$. Consider a local system~$\bbL$ of complex vector spaces on $U$ with finite monodromy group. For any point $u\in U(\bbC)$ with image $s=f(u)$ we denote by
\[
 G(u) \;=\; \mathrm{Im} \bigl(\pi_1(U_{s}, u) \; \longrightarrow \; \mathit{Aut}_\bbC(\bbL_u) \bigr)
\]
the monodromy group attached to the restriction $\bbL|_{U_s}$ of the local system $\bbL$ to the fibre $U_s = f^{-1}(s)$. Up to isomorphism this group only depends on $s=f(u)$.

\begin{lem} \label{lem:degenerate_monodromy}
Let $u_0 \in U(\bbC)$. Then there is a non-empty Zariski-open $S'\subseteq S$ with the following property: For all points $u\in f^{-1}(S'(\bbC))$, any identification $\bbL_{u_0} \cong \bbL_u$ of the stalks gives rise to an embedding of the monodromy groups such that the following diagram commutes.
\[
\xymatrix@C=2em@R=1.5em@M=0.6em{
 G(u_0) \ar@{^{(}-->}[r]^-\exists \ar@{^{(}->}[d] & G(u) \ar@{^{(}->}[d] \\
 \mathit{Aut}_\bbC(\bbL_{u_0}) \ar[r]^-\cong & \mathit{Aut}_\bbC(\bbL_u)
}
\]
\end{lem}

{\em Proof.} We may clearly replace the family $f: U\rightarrow S$ by its base change under any quasi-finite branched cover of a Zariski-open dense subset of $S$ containing $s_0$ if we also replace the local system $\bbL$ by its pull-back under this base change. Since \'etale-locally any smooth morphism admits a section, we can thus assume that our family admits a section $\sigma: S\rightarrow U$. By assumption $\bbL$ has finite monodromy, so after a further base change we can also assume that the local system $\sigma^*(\bbL)$ is trivial. Now recall~\cite[cor.~5.1]{Ve} that for any (not necessarily proper) morphism $f: U\rightarrow S$ of complex algebraic varieties there is a Zariski-open dense subset $S'\subseteq S$ over which the restriction $U' \rightarrow S'$ is a topologically locally trivial fibration in the analytic topology. For $s\in S'(\bbC)$ and $u=\sigma(s)$ we then have a split exact sequence 
\[
\xymatrix@C=2em@M=0.5em{
 \cdots \ar[r]
 & \pi_1(U_s, u) \ar[r]
 & \pi_1(U', u) \ar[r]
 & \pi_1(S, s) \ar[r] \ar@/_1.3pc/[l]|-{\, \sigma_* \,}
 & 1.
}
\]
Since $\sigma^*(\bbL)$ is trivial by construction, the lemma easily follows.
\qed

\medskip

\bigskip

{\em Acknowledgements.} This paper is based on chapter~5 of my Ph.D. thesis, and I would like to thank my advisor R.~Weissauer for his continuous support in all occuring questions and for many inspiring mathematical discussions.

\medskip

\bibliographystyle{amsplain}
\bibliography{BibliographyE6}

\end{document}